\documentclass[11pt]{article}
\usepackage{latexsym, graphicx, psfrag, xypic}
\usepackage{amscd, amsmath, amssymb, amsthm}
\usepackage{pstricks}

\newcommand{\cA}{\mathcal{A}}
\newcommand{\cB}{\mathcal{B}}
\newcommand{\cC}{\mathcal{C}}
\newcommand{\D}{\mathcal{D}}
\newcommand{\cE}{\mathcal{E}}
\newcommand{\cF}{\mathcal{F}}
\newcommand{\cG}{\mathcal{G}}
\newcommand{\caH}{\mathcal{H}}
\newcommand{\cK}{\mathcal{K}}
\newcommand{\cM}{\mathcal{M}}
\newcommand{\cN}{\mathcal{N}}
\newcommand{\LL}{\mathcal{L}}
\newcommand{\cO}{\mathcal{O}}
\newcommand{\cV}{\mathcal{V}}
\newcommand{\cEnd}{\mathcal{E}nd}
\newcommand{\cHom}{\mathcal{H}om}

\newcommand{\bC}{\mathbb{C}}

\newcommand{\bN}{\mathbb{N}}

\newcommand{\bP}{\mathbb{P}}

\newcommand{\bR}{\mathbb{R}}
\newcommand{\bZ}{\mathbb{Z}}

\newcommand{\RP}{\mathbb{RP}}

\newcommand{\ad}{\mathrm{ad}}

\newcommand{\Tr}{\mathrm{Tr}}
\newcommand{\Map}{\mathrm{Map}}
\newcommand{\Aut}{\mathrm{Aut}}
\newcommand{\End}{\mathrm{End}}
\newcommand{\rank}{\mathrm{rank}}

\newcommand{\diag}{\mathrm{diag}}

\newcommand{\vir}{\mathrm{vir}}

\newcommand{\fgl}{\mathfrak{gl}}
\newcommand{\fr}{\mathfrak{r}}

\newcommand{\fu}{\mathfrak{u}}
\newcommand{\fsu}{\mathfrak{su}}

\newcommand{\fm}{\mathfrak{m}}

\newcommand{\tSi}{ {\tilde{\Si}} }

\newcommand{\tP}{\tilde{P}}

\newcommand{\tg}{\tilde{g}}

\newcommand{\tbN}{\tilde{\bN}}

\newcommand{\ba}{\bar{a}}
\newcommand{\bb}{\bar{b}}
\newcommand{\bc}{\bar{c}}
\newcommand{\bd}{\bar{d}}

\newcommand{\bV}{\bar{V}}

\newcommand{\Si}{ {\Sigma} }
\newcommand{\ep}{\epsilon}

\newcommand{\ab}{a_1, b_1,\ldots, a_\ell, b_\ell}
\newcommand{\bab}{\ba_1,\bb_1,\ldots,\ba_\ell,\bb_\ell}
\newcommand{\pab}{\prod_{i=1}^\ell[a_i,b_i]}

\newcommand{\lra}{\longrightarrow}

\newcommand{\ymU}[2]{ X_{\mathrm{YM} }^{{#1},{#2}}(U(n)) }

\newcommand{\ZymU}[1]{ Z_{\mathrm{YM} }^{\ell, {#1}}(U(n)) }

\newcommand{\flS}[2]{ X_{\mathrm{flat}}^{ {#1}, {#2}}(U(1)) }
\newcommand{\ymS}[2]{ X_{\mathrm{YM} }^{{#1},{#2}}(U(1)) }

\newcommand{\ZymS}[1]{ Z_{\mathrm{YM} }^{\ell, {#1}}(U(1)) }
\newcommand{\kn}{{\frac{k}{n},\ldots,\frac{k}{n} }}
\newcommand{\km}{{\frac{k}{m},\ldots,\frac{k}{m} }}

\newcommand{\mkn}{ { -\frac{k}{n},\ldots,-\frac{k}{n} }}
\newcommand{\mtwokn}{ {-\frac{2k}{n},\ldots,-\frac{2k}{n} }}

\newcommand{\dbar}{\bar{\partial}}

\newcommand{\bG}{\mathcal{G}_0} 

\newcommand{\geqs}{\geqslant}
\newcommand{\heq}{\simeq}

\newcommand{\maps}{\longrightarrow}
\newcommand{\injects}{\hookrightarrow}

\newcommand{\isom}{\cong}
\newcommand{\cross}{\times}

\newcommand{\f}[1]{\underline{#1}}

\newcommand{\Css}{\mathcal{C}_{ss}}

\newcommand{\tr}{\mathrm{Tr}}
\newcommand{\hol}{\mathrm{hol}}

\newtheorem{thm}{Theorem}
\newtheorem{tm}[thm]{Theorem}
\newtheorem{df}[thm]{Definition}

\newtheorem{lm}[thm]{Lemma}
\newtheorem{rem}[thm]{Remark}

\newtheorem{pro}[thm]{Proposition}

\begin{document}

\title{Orientability in Yang-Mills Theory over Nonorientable Surfaces}
\author{\sc Nan-Kuo Ho, Chiu-Chu Melissa Liu, and Daniel Ramras}

%
%
%
%

\date{ }

\maketitle

\begin{abstract}
The first two authors have constructed a gauge-equivariant Morse stratification on the
space of connections on a principal $U(n)$-bundle
over a connected, closed, nonorientable surface $\Si$. This space can be identified with the real locus
of the space of connections on the pullback of this bundle over the orientable double cover
of $\Si$. In this context, the normal bundles to the Morse strata are real vector bundles.
We show that these bundles, and their associated homotopy orbit bundles, are orientable
for any $n$ when $\Si$ is not homeomorphic to the Klein bottle, and for
$n\leq 3$ when $\Si$ is the Klein bottle. We also derive similar orientability results when
the structure group is $SU(n)$.
\end{abstract}


\section{Introduction}\label{sec:introduction}

Consider a finite stratification $\{\cA_\mu\}$ of a manifold $S$.
If each stratum $\cA_\mu$ is a locally closed submanifold of $S$ with codimension $d_\mu$,
and the index set is partially ordered so that for any $\lambda$,
$$\bar{\cA}_\lambda\subset\bigcup_{\mu\geq\lambda}\cA_\mu$$
holds, then $\{\cA_\mu\}$ is called a Morse stratification.
A Morse stratification gives a Morse polynomial
$M_t(S;K)=\sum t^{d_\mu}P_t(\cA_\mu;K)$, where $P_t (-; K)$ denotes the Poincar\'{e} polynomial with coefficients in the field $K$.
The Morse inequalities state that there exists a polynomial $R_K(t)$ with
nonnegative coefficients such that $$M_t(S;K)=P_t(S;K)+(1+t)R_K(t).$$
Under fairly general conditions, the Morse inequalities hold for $K=\bZ_2$.
If, moreover, the normal bundle $\bN_\mu$ to each stratum $\cA_\mu$ is orientable,
 then these Morse inequalities hold for any coefficient field $K$.

Atiyah and Bott studied the moduli space of flat $G$-connections over
a Riemann surface via this Morse theoretical approach when the structure group $G$ is compact and connected.
One of their main results
is the computation of the $\cG$-equivariant Poincar\'{e} series
$P_t^\cG(\cA_{\mathrm{flat}};K)$ for the
space $\cA_{\mathrm{flat}}$ of flat connections on a principal bundle over a Riemann surface, where $\cG$ is the
gauge group.
They used the Yang-Mills functional, which is invariant
under the action of the gauge group, as a Morse-type function
and constructed a gauge equivariant
Morse stratification $\{\cA_\mu\}$ on the space $\cA$ of all connections
on a principal bundle over a Riemann surface. The
space $\cA_{\mathrm{flat}}$ of flat connections sits inside of the unique open
stratum $\cA_{ss}$ and is a deformation retract of
$\cA_{ss}$ via the Yang-Mills flow (Daskalopoulos \cite{Da}, R{\aa}de \cite{Rade}).
Thus, $\cA_{\mathrm{flat}}$ and $\cA_{ss}$ are homotopy equivalent,
and $P_t^\cG(\cA_{\mathrm{flat}};K)=P_t^\cG(\cA_{ss};K)$.
With this Morse stratification, one can write down the $\cG$-equivariant
Morse series of the space $\cA$ of all connections,
$$M_t^\cG(\cA;K)=\sum_{\mu\in I } t^{d_\mu}P_t^\cG(\cA_\mu;K),$$
and the $\cG$-equivariant Morse inequalities
$$M_t^\cG(\cA;K)=P_t^\cG(\cA;K)+(1+t)R_K(t),$$
where $d_\mu$ is the codimension of the stratum $\cA_\mu$, $I$ is
the index set of the stratification,
 and $R_K(t)$ is a power series with nonnegative coefficients.
In their construction,
the normal bundles $\bN_\mu$ are complex vector bundles, thus orientable, so
 $K$ can be any field.
In order to compute the $\cG$-equivariant Poincar\'{e} series
$P_t^\cG(\cA_{ss};K)$,
one needs four ingredients: $P_t^\cG(\cA;K)$, $d_\mu$,
$R_K(t)$, and $P_t^\cG(\cA_\mu;K)$ for all $\cA_\mu\neq \cA_{ss}$. Since the space $\cA$ of
all connections over a Riemann surface is an infinite dimensional
complex affine space and thus
contractible, $P_t^\cG(\cA;K)$ is just $P_t(B\cG;K)$, the Poincar\'{e} series
of the classifying space of $\cG$.
 As for $P_t^\cG(\cA_\mu)$, Atiyah and Bott found reduction formulas
\cite[Proposition 7.12]{ym} that reduce the question to smaller
groups.  The Morse index $d_\mu$ can be computed by Riemann-Roch
\cite[Equation (7.15)]{ym}.
Most importantly, they showed that this stratification
is $\cG$-equivariantly perfect \cite[Theorem 7.14]{ym}, i.e. $R_K(t)=0$, and
$$P_t(B\cG;K)=P_t^\cG(\cA;K)=\sum_{\mu\in I }
t^{d_\mu}P_t^\cG(\cA_\mu;K).$$
In the end, this method produces a recursive
formula for $P_t^\cG(\cA_{ss};K)$.

The first two authors defined a Yang-Mills functional
on the space of connections over any \emph{nonorientable} surface $\Si$ in \cite{HL1}.
Using this Yang-Mills functional, they constructed a $\cG$-equivariant Morse
stratification on the space of connections over $\Si$.
To be precise, consider the orientable double cover $\pi:\tSi\rightarrow\Si$,
and let $\tP=\pi^*P$ over $\tSi$ denote the pullback of
 a principal bundle $P$ over $\Si$. The non-trivial deck
 transformation of $\tSi$ induces an involution
on the space $\tilde{\cA}$ of connections of $\tP$ whose fixed point set
 is exactly the space $\cA$ of connections of $P$.
Ho and Liu define the Yang-Mills functional $L$ on $\cA$ to be
the restriction of the Yang-Mills functional $\tilde{L}$ on the
fixed point set of $\tilde{\cA}$. The absolute minimum of $L$ is
zero, achieved by flat connections on $P$. The gradient flow of
$L$ defines a $\cG$-equivariant Morse stratification $\{\cA_\mu\}$
on $\cA$. Indeed, the Morse stratification $\{\cA_\mu\}$ is just
the intersection of $\cA$ with the
Morse stratification $\{\tilde{\cA}_\mu\}$. This procedure also
tells us that the normal bundle $\bN_\mu$ to each stratum
$\cA_\mu$ in $\cA$  is the fixed locus of the normal bundle
$\tilde{\bN}_\mu$ to each stratum $\tilde{\cA}_\mu$ in $\tilde{\cA}$,
 which is complex (we will discuss in detail the various
involutions on vector bundles in Section \ref{sec:involution}). In
other words, the normal bundle $\bN_\mu$ to each Morse stratum
$\cA_\mu$ is a $\cG$-equivariant real vector bundle and hence is
not automatically orientable.

The $\cG$-equivariant Morse series of this stratification $\{\cA_\mu\}$ is
$$M_t^\cG(\cA;K)=\sum_{\mu\in I } t^{d_\mu}P_t^\cG(\cA_\mu;K).$$
Since the Yang-Mills strata admit gauge-invariant tubular neighborhoods
(see~\cite{Ra1} for a construction),
one can use the stratification to obtain the $\cG$-equivariant Morse inequalities
$$M_t^\cG(\cA;K)=P_t^\cG(\cA;K)+(1+t)R_K(t).$$
A priori, we cannot assume orientability of the normal bundles
$\bN_\mu$, so the Morse inequalities holds only for $K=\bZ_2$. To
compute the Poincar{\'e} series $P_t^\cG(\cA_{ss};K)$, we again
need four ingredients: $P_t^\cG(\cA;K)$, $d_\mu$, $R_K(t)$, and
$P_t^\cG(\cA_\mu;K)$ for all $\cA_\mu\neq\cA_{ss}$. Reduction
formulas for $P_t^\cG(\cA_\mu;K)$ and a formula for $d_\mu$ were given in
\cite{HL1, HL2}. On the
other hand, the computation of $P_t^\cG(\cA;K)$ is rather
difficult when $K=\bZ_2$ due to the existence of 2-torsion
elements in integral cohomology (see \cite[Section 5.3]{HL1}
\cite[Section 2]{HL2} for more details) and one is encouraged to
consider rational coefficients.  Hence we need to establish
orientability of the normal bundles.

Let $X_{h G}$ denote the
homotopy orbit space $EG\times_G X$. Then $\left(\bN_\mu\right)_{h
\cG}$ is also a real vector bundle over $\left(\cA_\mu\right)_{h \cG}$.
In this paper, we fix the structure group of the principal bundle $P$
to be the unitary group $U(n)$ or the special unitary group
$SU(n)$.  Our main result is:
\begin{thm}\label{thm:main}
Suppose that either
\textup{(i)} $\chi(\Si) = 0$ (so that $\Si$ is homeomorphic to the Klein bottle)
and $n\leq 3$, or
\textup{(ii)} $\chi(\Si)\neq 0$ and $n$ is any positive integer.
Then
$\left(\bN_\mu\right)_{h \cG}$ is an orientable vector
bundle over $\left(\cA_\mu\right)_{h \cG}$ for all $\mu$.
As a consequence,  $\bN_\mu$ is an orientable
vector bundle over $\cA_\mu$ for all $\mu$.
\end{thm}

In \cite{HL3}, the first two authors discuss how far
this stratification is from being perfect, i.e. what the power series $R_K(t)$ looks like.
They define the notion of {\em antiperfection}, which leads to
some conjectural formulas for $P_t^\cG(\cA_{ss};K)$.

Thomas Baird \cite{B} has recently proven the
 formula conjectured in \cite{HL3} for the $\cG$-equivariant Poincar\'{e} series of
the space of flat $U(3)$-connections over a non-orientable surface.
His argument relies on Yang-Mills theory, and in particular
uses our orientability results.
Thus Baird's work may be viewed as a concrete application of the results in this paper.

\section{Preliminaries}\label{sec:prelim}

Let $\tSi$ be a Riemann surface.
Let $P^{n,k}_\tSi$ denote the degree $k$ principal $U(n)$-bundle
on $\tSi$. Let $\rho: U(n) \to  GL(n,\bC)$ be the fundamental representation,
and let $E= P^{n,k}_\tSi \times_\rho \bC^n$ be the associated
complex vector bundle over $\tSi$. Then $E$ is a rank $n$, degree $k$
complex vector bundle equipped with a Hermitian metric $h$, and
the unitary frame bundle $U(E,h)$ of the Hermitian
vector bundle $(E, h)$ is isomorphic to $P^{n,k}_\tSi$
as a $C^\infty$ principal $U(n)$-bundle.

\subsection{Hermitian,  $(0,1)$-, and $(1,0)$-connections}
\label{sec:hermitian}

Let $\cA(P^{n,k}_\tSi)$ be the space of $U(n)$-connections
on $P^{n,k}_\tSi$, which can be identified with $\cA(E,h)$,
the space of Hermitian connections on $(E,h)$ (connections
on $E$ which are compatible with the Hermitian structure $h$,
cf. \cite[pp.76]{We}) .
It is a complex affine space whose vector
space of translations is $\Omega^1_\tSi(\ad P^{n,k}_\tSi)$, where
the complex structure is given by the Hodge star $*$ (cf. \cite{ym}).
Let $\cC(E)$ denote the space of $(0,1)$-connections
$\dbar: \Omega^0_{\tSi}(E)\to \Omega^{0,1}_\tSi(E)$,
and let $\cC'(E)$ denote the space of $(1,0)$-connections
$\partial: \Omega^0_{\tSi}(E)\to \Omega^{1,0}_\tSi(E)$.
Recall that a $(0,1)$-connection  (resp. $(1,0)$-connection)
defines a holomorphic (resp. anti-holomorphic) structure
on $E$ if and only if $\dbar^2 =0$ (resp. $\partial^2 =0$)
(cf. \cite[Section 2.2.2]{DK}); now
$\Omega_\tSi ^{0,2} = 0$ (resp. $\Omega_\tSi^{2,0}=0$)
since $\dim_\bC\tSi=1$, so the
integrability condition $\dbar^2=0$ (resp. $\partial^2=0$)
holds automatically.
The local holomorphic (resp. anti-holomorphic) sections
are solutions to $\dbar s=0$ (resp. $\partial s=0$).
$\cC(E)$ and $\cC'(E)$ are complex
affine spaces whose vector spaces of translations
are $\Omega^{0,1}_{\tSi}(\End(E))$ and
$\Omega^{1,0}_{\tSi}(\End(E))$, respectively (cf. \cite{ym}).

Given a Hermitian connection $\nabla: \Omega^0_{\tSi}(E)\to \Omega^1_{\tSi}(E)$,
let $\nabla': \Omega^0_\tSi(E)\to \Omega^{1,0}_{\tSi}(E)$ and
$\nabla'':\Omega^0_\tSi(E)\to \Omega^{0,1}_{\tSi}(E) $ be the $(1,0)$ and
$(0,1)$ parts of $\nabla$.  Then $\nabla\mapsto \nabla''$
and $\nabla\mapsto \nabla'$ define
isomorphisms $j: \cA(P^{n,k}_{\tSi})\to \cC(E)$ and
$j':\cA(P^{n,k}_{\tSi})\to \cC'(E)$ of real affine spaces.
Their differentials
$$
j_*: \Omega^1_\tSi(\ad P^{n,k}_\tSi)\to \Omega^{0,1}_\tSi(\End E),\quad
j'_*: \Omega^1_\tSi(\ad P^{n,k}_\tSi)\to \Omega^{1,0}_\tSi(\End E),
$$
are complex linear and conjugate linear, respectively.
More explicitly, $j_*$ and $j'_*$ are $C^\infty(\tSi,\bR)$-linear,
so they are induced by real vector bundle maps
$\tilde{j}: T^*_{\tSi}\otimes \ad P^{n,k}_{\tSi}\to (T^*_{\tSi})^{0,1}\otimes \End E$
and $\tilde{j}':T^*_{\tSi}\otimes \ad P^{n,k}_{\tSi}\to (T^*_{\tSi})^{1,0}\otimes \End E$,
respectively. Given any point $x\in \tSi$, let $dz=dx + idy$ be a local basis
of $(T^*_{\tSi})_p^{1,0}$ and let $X,Y\in \fu(n)$.  Then $X+ iY \in \fgl(n,\bC)$, and
$$
\tilde{j}(X dx + Y dy) = \frac{1}{2}(X+iY) d\bar{z},\quad
\tilde{j}'(X dx + Y dy) =\frac{1}{2} (X-iY) d z.
$$
The complex structure on $(T^*_\tSi\otimes \ad P^{n,k}_\tSi)_x$
is given by the Hodge star: $*(X dx + Y dy) = - Y dx + X dy$. It is
straightforward to check that $\tilde{j}$ is complex
linear and $\tilde{j}'$ is conjugate linear.

Given a $(0,1)$-connection $\dbar$ on a Hermitian vector bundle $(E,h)$ over
$\tSi$, there is a unique connection $\nabla$ on $E$ which is compatible with
$h$ and such that $\nabla''=\dbar$ (see e.g. \cite[pp.78]{We}). We denote
this canonical Hermitian connection by $\nabla_{h,\dbar}$.
The map $j'\circ j^{-1}:\cC(E)\to \cC'(E)$ is given
by $\dbar\mapsto (\nabla_{h,\dbar})'$, where $(\nabla_{h,\dbar})'$ is the
$(1,0)$-part of $\nabla_{h,\dbar}$.

Let $E^\vee$ be the complex dual of $E$ (see e.g.
\cite[pp.168-169]{MS}). Then $E^\vee$ is a rank $n$, degree $-k$
complex vector bundle equipped with a Hermitian metric $h^\vee$
induced by $h$. More explicitly, if $\{ e_1,\ldots,e_n\}$ is a local
orthonormal frame of the Hermitian vector bundle $(E,h)$, then its
dual coframe $\{ e_1^\vee, \ldots, e_n^\vee\}$ is a local
orthonormal frame of the Hermitian vector bundle $(E^\vee, h^\vee)$.
The map $v \mapsto h(\cdot, v)$ defines a conjugate linear bundle
map $E\to E^\vee$ which induces an isomorphism $I_h: (E,h)\cong
(\overline{E^\vee}, \overline{h^\vee})$ of Hermitian vector bundles.
We have $U(\overline{E^\vee}, \overline{h^\vee})\cong U(E,h) \cong
P^{n,k}_\tSi$.

A $(0,1)$-connection $\dbar$ on $E$ induces a $(0,1)$-connection
$\dbar^\vee$ on $E^\vee$ and a $(1,0)$-connection $\partial^\vee$ on
$\overline{E^\vee}$. This gives a map $j_1:\cC(E)\to \cC'(\overline{E^\vee})$.
The map $v\mapsto h(\cdot, v)$ defines an isomorphism
$E\cong \overline{E^\vee}$ of $C^\infty$ complex vector bundles,
which induces an isomorphism $j_2: \cC'(E)\to \cC'(\overline{E^\vee})$
of complex affine spaces.
It is straightforward to check that
$j'\circ j^{-1} =j_2^{-1}\circ j_1:\cC(E)\to \cC'(E)$.

\subsection{Gauge groups}\label{sec:gauge}

Let $GL(E)$ be the frame bundle of the complex vector bundle $E$.
Let $U(E,h)$ be the unitary frame bundle of the Hermitian vector
bundle $(E,h)$ as in the previous subsection.
Then $GL(E)$ is a principal $GL(n,\bC)$-bundle over $\tSi$, and
$U(E,h)$ is a principal $U(n)$-bundle over $\tSi$.
Let $\Aut(E)$ be the (infinite dimensional) group of complex vector bundle isomorphisms $E\to E$,
and let $\Aut(E,h)$ be the (infinite dimensional) group of Hermitian bundle isomorphisms $(E,h)\to (E,h)$. (See \cite[Section 2]{ym} for details.)
Then  $\Aut (E)\cong \Aut GL(E)$ and $\Aut(E,h)\cong \Aut U(E,h)$;
$\Aut(E,h)$ is a subgroup of $\Aut (E)$.

$\Aut(E)$ acts on $\cC(E)$ by  $u\cdot \dbar = u\circ \dbar \circ u^{-1}$
and $\Aut(E,h)$ acts on $\cA(E,h)$ by $u\cdot \nabla = u\circ \nabla\circ u^{-1}$.
More explicitly, relative to a local orthonormal frame,
a $(0,1)$-connection on $E$ is of the form
$$
\dbar= \dbar_0 + B,
$$
where $\dbar_0$ is the usual Cauchy-Riemann operator and $B$ is
a $\fgl(n,\bC)$-valued $(0,1)$-form; a unitary connection is of the
form
$$
\nabla = d + A,
$$
where $d$ is the usual exterior derivative and $A$ is a
$\fu(n)$-valued $1$-form. An element $u$ in the gauge group
$\Aut(E)$ is locally a $GL(n,\bC)$-valued function, and acts on the form
$B$ by
\begin{equation}\label{eqn:uB}
B\mapsto  u B u^{-1} -(\dbar_0 u)  u^{-1};
\end{equation}
an element $u$ in the gauge group $\Aut(E,h)$ is locally a
$U(n)$-valued function, and acts on the form $A$ by
\begin{equation}\label{eqn:uA}
A\mapsto  u A u^{-1} -(du) u^{-1}.
\end{equation}
In particular, if $u\in GL(n,\bC)$ (resp.
$U(n)$) is a {\em constant} gauge transformation, then it acts on $B$
(resp. $A$) by $B\mapsto  u B u^{-1}$
(resp. $A\mapsto  u A u^{-1}$ ).

Given $u\in \Aut(E)$ and $x\in \tSi$, $u_x: E_x\to E_x$ is a complex
linear isomorphism for all $x\in \tSi$. The dual of $u_x$ is a
complex linear isomorphism $(u_x)^\vee: (E_x)^\vee\to (E_x)^\vee
=(E^\vee)_x$. It induces a complex linear isomorphism
$$
\overline{(u_x)^\vee}: \overline{(E_x)^\vee}\to \overline{(E_x)^\vee}\cong (\overline{E^\vee})_x.
$$
Define $\overline{u^\vee}\in \Aut(\overline{E^\vee})$ by $(\overline{u^\vee})_x
=\overline{(u_x)^\vee}$. Then $u\maps \overline{u}^\vee$ defines a group homomorphism
$\Aut(E)\to \Aut(\overline{E^\vee})$. The isomorphism $I_h: E\cong \overline{E^\vee}$ allows us to identify
$\Aut(E)$ with $\Aut(\overline{E^\vee})$. We let $\tilde{I}_h:\Aut(E)\to \Aut(\overline{E^\vee})$
be this $h$-dependent identification, and let
$\phi_h:\Aut(E)\to \Aut(E)$ be defined by $u\mapsto \tilde{I}_h(\overline{u^\vee})$.
Then $\phi_h$ can be described explicitly as follows. Let $u\in \Aut(E)$, and let
$A\in GL(n,\bC)$ be the matrix of $u_x:E_x\to E_x$ with respect to
an orthonormal basis of $(E_x, h_x)$.
Then $\phi_h(u)_x = \overline{(A^t)^{-1}}$.
Note that $\phi_h:\Aut(E)\to \Aut(E)$ is an involution, and the
fixed locus $\Aut(E)^{\phi_h}=\Aut(E,h)$.

\section{Involution} \label{sec:involution}
Let $\Si$ be a closed nonorientable surface, and let $\pi:\tSi\to \Si$ be its orientable double
cover.   Then $\tSi$ is a Riemann surface, and the non-trivial deck transformation is an anti-holomorphic,
anti-symplectic involution $\tau:\tSi\to \tSi$ such that $\pi\circ\tau =\pi$.

\subsection{The action of $\tau$ on holomorphic structures} \label{sec:tau-cC}
There is an anti-holomorphic,
anti-symplectic map $\tau_{\cA}:\cA(E,h) \to \cA(\tau^*E,\tau^*h)$
given by $\nabla\mapsto \tau^*\nabla$. Note that
$\tau^* P_{\tSi}^{n,k}\cong P_{\tSi}^{n,-k}$, so
$\cA(\tau^*E, \tau^*h)\cong \cA(P^{n,-k}_\tSi)$.
We have
$$
(\tau^*\nabla)'= \tau^*(\nabla''),\quad (\tau^*\nabla)'' =\tau^*(\nabla'),
$$
so there are maps
\begin{eqnarray*}
&& \tau^* : \cC(E)\to \cC'(\tau^* E), \quad \cC(\tau^*E)\to \cC'(E),\quad \dbar \mapsto \tau^* \dbar\\
&& \tau^*:  \cC'(E)\to \cC(\tau^* E), \quad \cC'(\tau^*E)\to \cC(E),\quad \partial \mapsto \tau^* \partial,
\end{eqnarray*}
such that $\tau^* \circ \tau^*$ is the identity map.

Define $\tau_\cC := j\circ \tau_{\cA} \circ j^{-1} =\tau^* \circ j'\circ j^{-1}: \cC(E)\to \cC(\tau^*E)$.
Then $\tau_\cC$ is given by $\dbar \mapsto \tau^* (\nabla_{h,\dbar})'$. In the rest of this subsection, we
study the effect of $\tau_\cC$ on the Harder-Narasimhan filtration.

Let $\cE$ denote $E$ equipped with a $(0,1)$-connection (holomorphic structure),
so that $\cE$ can be viewed as a point in $\cC(E)$. Let
$$
0=\cE_0\subset \cE_1 \subset \cdots\subset  \cE_r=\cE
$$
be the Harder-Narasimhan filtration, so
$\cE_j/\cE_{j-1}$ is semi-stable. Set
$$
\D_j =\cE_j/\cE_{j-1},\quad n_j=\rank_\bC \D_j,\quad
k_j=\deg \D_j.
$$
The Atiyah-Bott type of $\cE$ is
\begin{equation}\label{eqn:ABtype}
\mu(\cE) =\Bigl(\underbrace{\frac{k_1}{n_1},\ldots,\frac{k_1}{n_1}}_{n_1},\ldots,
\underbrace{\frac{k_r}{n_r},\ldots,\frac{k_r}{n_r}}_{n_r}\Bigr),\quad
\textup{where }\frac{k_1}{n_1} >\cdots > \frac{k_r}{n_r}.
\end{equation}

Recall that $\rank_\bC E= n$ and $\deg E=k$. Let
$$
I_{n,k}=\Bigl\{
\Bigl(\underbrace{\frac{k_1}{n_1},\ldots,\frac{k_1}{n_1}}_{n_1},\ldots,
\underbrace{\frac{k_r}{n_r},\ldots, \frac{k_r}{n_r}}_{n_r}\Bigr)\
\Bigr|\ \frac{k_1}{n_1}>\cdots > \frac{k_r}{n_r},\  \sum_{j=1}^{r}
n_j =n,\ \sum_{j=1}^r k_j =k \Bigr\},
$$
and for $\mu\in I_{n,k}$ let $\cC_\mu = \{\cE \in \cC(E) \ \Bigr|\ \mu(\cE) = \mu\}$.
The Harder-Narasimhan strata of $\cC(E)$ are $\{ \cC_\mu\mid \mu\in I_{n,k}\}$ (cf. \cite[Section 7]{ym}).

Using the isomorphism $I_h: E\cong \overline{E^\vee}$ defined in Section \ref{sec:hermitian},
we may identify $\cC(\tau^*E)$ with $\cC(\tau^*\overline{E^\vee})$. Then
$\tau_\cC:\cC(E)\to \cC(\tau^*\overline{E^\vee})$ is given by
$\cE \mapsto \tau^* \overline{\cE^\vee}$,
where $\cE, \cE^\vee, \tau^*\overline{\cE^\vee}$ are holomorphic vector bundles
over $\tSi$, while
$\overline{\cE^\vee}$ is an anti-holomorphic vector bundle over $\tSi$.

For $j=0,\ldots,r$, define a holomorphic subbundle
$(\cE^\vee)_{-j}$ of $\cE^\vee$ by
$$
\bigl((\cE^\vee)_{-j}\bigr)_x =\{ \alpha\in \cE^\vee_x \mid \alpha(v)=0\ \forall v\in (\cE_j)_x \}.
$$
Then $(\cE^\vee)_{-j}=(\cE/\cE_j)^\vee$. The Harder-Narasimhan
filtration of $\cE^\vee\in \cC(E^\vee)$ is given by
$$
0=(\cE^\vee)_{-r}\subset (\cE^\vee)_{-(r-1)}\subset \cdots
\subset (\cE^\vee)_{-1}\subset (\cE^\vee)_0=\cE^\vee
$$
Notice that
$$(\cE^\vee)_{-i}/(\cE^\vee)_{-(i+1)} \cong (\cE_{i+1}/\cE_i)^\vee=(\D_{i+1})^\vee.$$
Set $\caH_j= (\cE^\vee)_{-(r-j)}/(\cE^\vee)_{-(r-j+1)}$. Then
$$
\caH_j \cong (\D_{r+1-j})^\vee,\quad
\rank_\bC \caH_j= n_{r+1-j},\quad \deg\caH_j = -k_{r+1-j}.
$$
Hence the Atiyah-Bott type of $\cE^\vee$ is
$$
\mu=\Bigl(\underbrace{-\frac{k_r}{n_r},\ldots,-\frac{k_r}{n_r}}_{n_r},\ldots,
\underbrace{-\frac{k_1}{n_1},\ldots,-\frac{k_1}{n_1}}_{n_1}\Bigr),\quad
\textup{where }
-\frac{k_r}{n_r} >\cdots > -\frac{k_1}{n_1}.
$$

For $j=0,\ldots,r$, define a holomorphic subbundle
$\tau_\cC(\cE)_{-j}$ of $\tau_\cC(\cE)=\tau^*\overline{\cE^\vee}$ by
$\tau_\cC(\cE)_{-j} = \tau^*\overline{(\cE^\vee)_{-j}}$. The
Harder-Narasimhan filtration of $\tau_\cC(\cE)$ is given by
$$
0=\tau_\cC(\cE)_{-r}\subset
\tau_\cC(\cE)_{-(r-1)}\subset \cdots
\subset \tau_\cC(\cE)_0=\tau_\cC(\cE)
$$
Let $\cK_j =\tau_\cC(\cE)_{-(r-j)}/\tau_\cC(\cE)_{-(r-j+1)}$.
Then
$$
\cK_j \cong \tau^*\overline{\caH_j} \cong \tau^*\overline{(\D_{r+1-j})^\vee}=\tau_\cC(\D_{r+1-j}),
$$
$$\textrm{and} \quad \rank_\bC \cK_j = n_{r+1-j},\quad \deg\cK_j = -k_{r+1-j}.
$$
The Atiyah-Bott type of $\tau_\cC(\cE)$ is
$$
\mu=\Bigl(\underbrace{-\frac{k_r}{n_r},\ldots,-\frac{k_r}{n_r}}_{n_r},\ldots,
\underbrace{-\frac{k_1}{n_1},\ldots,-\frac{k_1}{n_1}}_{n_1}\Bigr).
$$

From the above discussion, we conclude:
\begin{lm}
Let $\tau_\cC:\cC(E)\to\cC(\tau^*E)\cong \cC(\tau^*\overline{E^\vee})$
be defined as above, and define
$\tau_\cC: \Aut (E) \to \Aut(\tau^*E)$ by $u\mapsto \tau^* \phi_h(u)$, where
$\phi_h$ is defined as in Section \ref{sec:gauge}.
Define $\tau_0:I_{n,k}\to I_{n,-k}$ by
$$
\Bigl(\underbrace{\frac{k_1}{n_1},\ldots,\frac{k_1}{n_1}}_{n_1},\ldots,
\underbrace{\frac{k_r}{n_r},\ldots, \frac{k_r}{n_r}}_{n_r}\Bigr)
\mapsto
 \Bigl(\underbrace{-\frac{k_r}{n_r},\ldots,-\frac{k_r}{n_r}}_{n_r},\ldots,
\underbrace{-\frac{k_1}{n_1},\ldots, -\frac{k_1}{n_1}}_{n_1}\Bigr),
$$
Then
\begin{enumerate}
\item $\tau_\cC:\cC(E)\to \cC(\tau^*E)\cong \cC(\tau^*\overline{E^\vee})$ maps $\cC_\mu$ bijectively
to $\cC_{\tau_0(\mu)}$.
\item $\tau_\cC$ is equivariant with
respect to the $\Aut(E)$-action on $\cC(E)$ and $\Aut(\tau^*E)$-action on $\cC(\tau^*E)$, i.e.,
$$
\tau_\cC(u\cdot \dbar)= \tau_\cC(u)\cdot \tau_\cC(\dbar),\quad
u\in \Aut(E),\quad \dbar\in \cC(E).
$$
\end{enumerate}
\end{lm}
\subsection{The degree zero case}

Let $P\to \Si$ be a principal $U(n)$-bundle, and let $\tP =\pi^*
P$ be the pull back principal $U(n)$-bundle on $\tSi$. We first
review some facts about $\tP$ (see \cite[Section 3.2]{HL1} for
details). The pull back bundle $\tP\cong P^{n,0}_\tSi \cong
\tSi\times U(n)$ is topologically trivial. We wish to describe an
involution $\tilde{\tau}: \tP\to \tP$ which is $U(n)$-equivariant,
covers the involution $\tau:\tSi\to\tSi$, and satisfies $P
=\tP/\tilde{\tau}_s$. Fixing a trivialization $\tP \cong
\tSi\times U(n)$, any such involution must be given by
$\tilde{\tau}_s: \tSi\times U(n)\to \tSi\times U(n)$,
$(x,h)\mapsto (\tau(x), s(x)h)$, for some  $C^\infty$ map
$s:\tSi\to U(n)$ satisfying $s(\tau(x))=s(x)^{-1}$.

The topological type of a principal $U(n)$-bundle $P\to \Si$ is
classified by  $c_1(P) \in H^2(\Si;\bZ)\cong \bZ/2\bZ$. Let
$P^{n,+}_\Si$ and $P^{n,-}_\Si$ denote the principal $U(n)$-bundles
on $\Si$ with $c_1=0$ and $c_1=1$ in $\bZ/2\bZ$, respectively. Let
$\tau_\ep$ be the involution on $P^{n,0}_\tSi =\tSi\times U(n)$
defined by a constant map $s(x)=\ep \in U(n)$.  We must have
$\ep^2 =I_n$, so $\det\ep=\pm 1$. Then $P^{n,0}_\tSi/\tau_\ep
\cong P^{n,\pm}_\Si$ if $\det\ep =\pm 1$. We choose $\ep_\pm$ to
be the diagonal matrix $\diag(\pm 1, 1, \ldots, 1)$, and define
$\tau^\pm = \tau_{\ep_\pm}$. Then $P^{n,0}_\tSi/\tau^\pm \cong
P^{n,\pm}_\Si$.

Let $E=P^{n,0}_\tSi \times_\rho \bC^n \cong \tSi\times \bC^n$,
where $\rho:U(n)\to GL(n,\bC)$ is the fundamental representation.
Then $\tau^\pm$ induces
an involution $\tau^\pm: E\cong \tSi\times \bC^n \to E \cong \tSi\times \bC^n$
given by $(x,v)\mapsto (\tau(x), \ep_\pm v)$.  The two involutions $\tau^+, \tau^-$
give two isomorphisms
$\tau^* E\cong E$, which induce isomorphisms
$$
\cA(E,h) \cong \cA(\tau^*E,\tau^*h),\quad
\cC(E)\cong \cC(\tau^*E),\quad
\Aut(E)\cong \Aut(\tau^*E).
$$
Therefore, we have involutions
$$
\tau^\pm_\cA: \cA(E,h)\rightarrow \cA(E,h),\quad \tau^\pm_\cC:
\cC(E)\rightarrow \cC(E),\quad
\tau^\pm_\cC:\Aut(E)\rightarrow \Aut(E),
$$
and $\tau^\pm_\cC:\cC(E)\to \cC(E)$ is $\Aut(E)$-equivariant
with respect to the $\Aut(E)$-action on $\cC(E)$.
We have
$$
\cA(P^{n,\pm}_\Si) =\cA(E,h)^{\tau_\cA^\pm} \cong \cC(E)^{\tau_\cC^\pm},\quad
\Aut(P^{n,\pm}_\Si) \cong \Aut(E,h)^{\tau_\cC^\pm},
$$
where $\Aut(E,h)\subset \Aut(E)$ is the group of unitary gauge
transformations of the Hermitian vector bundle $(E,h)$.  The
following two equivariant pairs are isomorphic:
$$
\left(\cA(P^{n,\pm}_\Si), \Aut(P^{n,\pm}_\Si)\right)
\cong \left(\cC(E)^{\tau_\cC^\pm},\Aut(E,h)^{\tau_\cC^\pm}\right).
$$

\subsection{$SU(n)$-connections}
Let $Q^n_\Si \to \Si$ be a principal $SU(n)$-bundle. Then
$Q^n_\Si $ is topologically trivial.
We fix a trivialization $Q^n _\Si \cong \Si\times SU(n)$, which
allows us to identify the space $\cA(Q^n_\Si)$ of connections
on $Q^n_\Si$ with the vector space of $\fsu(n)$-valued 1-forms
on $\Si$. Let $P^{n,+}_\Si \cong \Si\times U(n)$ be the trivial $U(n)$-bundle
on $\Si$, as before. The short exact sequence of vector spaces
$$
0\to \fsu(n)\to \fu(n)\stackrel{\tr}{\to} \fu(1)\to 1
$$
induces a short exact sequence of infinite dimensional  vector spaces
$$
0\to \cA(Q^n_\Si) \to \cA(P^{n,+}_\Si) \stackrel{\tr}{\to} \cA(P^{1,+}_\Si)\to 0.
$$
The Yang-Mills functional on $\cA(Q^n_\Si)$ is the restriction of the Yang-Mills functional
on $\cA(P^{n,+}_\Si)$.  The Morse stratifications on
$\cA(P^{n,+}_\Si)$ and on $\cA(Q^n_\Si)$ are given by
$$
\cA(P^{n,+}_\Si)=\bigcup_{\mu\in I}\cA_\mu,\quad
\cA(Q^n_\Si) =\bigcup_{\mu\in I} \cA'_\mu
$$
where $\cA'_\mu = \cA_\mu\cap \cA(Q^n_\Si)$ is nonempty
for any $\mu\in I$. Given $\mu\in I$ such that
$\cA_\mu \neq \cA_{ss}$, let $\bN_\mu$ (resp. $\bN_\mu'$) be the normal bundle of
$\cA_\mu$ (resp. $\cA_\mu'$) in $\cA(P^{n,+}_\Si)$ (resp.
$\cA(Q^n_\Si)$). Let $\iota_\mu:\cA'_\mu\hookrightarrow \cA_\mu$
be the inclusion. Then $\bN'_\mu =\iota_\mu^* \bN_\mu$.
Therefore, if $\bN_\mu$ is an orientable real vector bundle over $\cA_\mu$
then $\bN_\mu'$ is an orientable real vector bundle over $\cA_\mu'$.

The short exact sequence of Lie groups
$$
1\to SU(n)\to U(n)\stackrel{\det}{\to} U(1)\to 1
$$
induces a short exact sequence of infinite dimensional gauge groups
$$
1\to \Map(\Si,SU(n)) \to \Map(\Si,U(n)) \to \Map(\Si, U(1)) \to 1
$$
or equivalently,
\begin{equation}\label{eqn:USU-cG}
1\to \Aut(Q^n_\Si)\to \Aut(P^{n,+}_\Si) \to \Aut(P^{1,+}_\Si)\to 1.
\end{equation}
In particular, $\cG':= \Aut(Q^n_\Si)$ is a subgroup
of $\cG :=\Aut(P^{n,+}_\Si)$; indeed $\cG'$ is a closed, normal subgroup of $\cG$.
We have the following diagram:
$$
\begin{CD}
E\cG\times_{\cG'} \bN_\mu' @>>> E\cG\times_{\cG'} \bN_\mu @>>> E\cG\times_\cG \bN_\mu\\
@VVV @VVV @VVV\\
E\cG\times_{\cG'} \cA_\mu' @>>> E\cG\times_{\cG'} \cA_\mu  @>>> E\cG\times_\cG \cA_\mu
\end{CD}
$$
which can be identified with
$$
\begin{CD}
(\bN_\mu')_{h\cG'} @>>> (\bN_\mu)_{h\cG'} @>>> (\bN_\mu)_{h\cG}\\
@VVV @VVV @VVV\\
 (\cA_\mu')_{h\cG'} @>{\iota_\mu }>> (\cA_\mu)_{h\cG'}  @>{q_\mu }>> (\cA_\mu)_{h\cG}
\end{CD}
$$
where
$$
(\bN_\mu)_{h\cG'}= q_\mu^* \left((\bN_\mu)_{h\cG}\right),\quad
(\bN_\mu')_{h\cG'} = \iota_\mu^*\left((\bN_\mu)_{h\cG'}\right).
$$
Therefore, if $(\bN_\mu)_{h\cG}$ is an orientable vector bundle
over $(\cA_\mu)_{h\cG}$ then $(\bN_\mu')_{h\cG'}$ is an orientable
vector bundle over $(\cA_\mu')_{h\cG'}$.

From the above discussion, if Theorem \ref{thm:main}
holds for $U(n)$ then it holds for $SU(n)$.
In the remainder of this paper, we prove
Theorem \ref{thm:main}  for $U(n)$.

\section{Reduction} \label{sec:reduction}

Let $\Si$ denote a closed, non-orientable surface.
In this section, we reduce the question of orientability
for normal bundles of Morse strata in $\cA(P^{n, \pm}_\Si)$
to the question of orientability for certain real vector bundles
$V_{n,k}$ over the representation varieties
associated to central Yang-Mills connections on
$P^{n,k}_\tSi$.

The reduction will pass through a variety of gauge-theoretical
spaces, most of which are not CW complexes.  Hence one needs to be
careful in applying the usual bundle-theoretical arguments.  In the
end, however, we will show that the normal bundle to each Morse
stratum $\cA_\mu$, when considered equivariantly as a bundle over
$\left(\cA_\mu\right)_{h\cG}$, is pulled back under a weak
equivalence from a bundle over the homotopy orbit space
$\left(\cN_\mu/\bG\right)_{hU(n)}$ (here $\cN_\mu$ denotes the set
of type $\mu$ Yang-Mills connections).  These representation
varieties are analytical sets~\cite{HL1}, and their homotopy orbit
spaces are triangulable by results of Illman~\cite{Il}. The fact
that the normal bundle is pulled back from a bundle over a CW
complex will allow us to use standard bundle-theoretical arguments.
At the end of this section we will summarize the arguments to follow,
so as to make the overall strategy of the reduction clear and the
proof rigorous.

\subsection{Reduction to Levi subgroups}
\label{sec:levi}

On the stratum $\cC_\mu$,
where $\mu$ is as in Equation \eqref{eqn:ABtype}, we
will proceed to reduce the $U(n)$-gauge group to a
Levi subgroup corresponding to $U(n_1)\times \cdots \times U(n_r)$.  Our arguments follow \cite[Section 7]{ym} closely.

Let $\cF_\mu$ denote the
space of all $C^\infty$ filtrations of type $\mu$. The Harder-Narasimhan filtration provides a continuous map
$p: \cC_\mu\to \cF_\mu$. Let $E_\mu\in \cF_\mu$ be a fixed $C^\infty$ filtration
of $E$ and let $\cB_\mu = p^{-1}(E_\mu)$. We choose splittings of the
filtration $E_\mu$ to obtain a direct sum decomposition $E^0_\mu=D_1\oplus \cdots\oplus D_r$
of $E$, and let $\cB_\mu^0\subset \cB_\mu$ be the space of complex structures
compatible with the direct sum decomposition $E_\mu^0$.  The inclusion $\cB_\mu^0 \injects \cB_\mu$
splits the fibration $\cB_\mu \to \cB_\mu^0$, which has a vector space as fiber.  Hence this inclusion is a weak equivalence.
Since $\cC_\mu$ is the extension of the $\Aut(E_\mu)$-space $\cB_\mu$ to a $\Aut(E)$-space, we have a
homeomorphism of homotopy orbit spaces $\left(\cC_\mu\right)_{ h \Aut(E)} \isom \left(\cB_\mu\right)_{h\Aut(E_\mu)}$.  Thus
$$
\left(\cC_\mu\right)_{h \Aut(E)} \isom \left(\cB_\mu\right)_{h\Aut(E_\mu)} \sim (\cB_\mu^0)_{h \Aut(E_\mu^0)}
\isom \prod_{j=1}^r \cC_{ss}(D_i)_{h \Aut(D_i)}.
$$

Let $\tbN_\mu\to \cC_\mu$ be the normal bundle of $\cC_\mu$ in $\cC(E)$.
Given $\cE\in \cB^0_\mu \subset \cC_\mu$, $\cE$ is a direct sum of holomorphic
subbundles $\D_1, \ldots, \D_r$, and
$$
(\tbN_\mu)_\cE = \bigoplus_{i<j} H^1(\tSi, \cHom(\D_i, \D_j)).
$$

The Harder-Narasimhan filtration is again a continuous map $p:\cC_{\tau_0(\mu)}\to \cF_{\tau_0(\mu)}$. We have
$\tau^* E_\mu \in \cF_{\tau_0(\mu)}$. Let $\cB_{\tau_0(\mu)} = p^{-1} (\tau^* E_\mu)$
and let $\cB^0_{\tau_0(\mu)}$ be the space of complex structures
compatible with the direct sum decomposition $\tau^*(E^0_\mu)=\tau^* D_1\oplus \cdots\oplus \tau^* D_r$.
Then $\tau(\cB_\mu^0)=\cB_{\tau_0(\mu)}^0$.

Given a holomorphic vector bundle $\cV\to \tSi$,
let $\cO(\cV)$ be the sheaf of local holomorphic
sections on $\cV$. Then for $i=0,1$,
\begin{equation}\label{eqn:cohomology}
H^i_{\dbar}(\tSi,\cV) \cong H^i(\tSi,\cO(\cV)) \cong \check{H}^i(\{ U_\alpha\}, \cO(\cV))
\end{equation}
where $H^i_{\dbar}(\tSi,\cV)$ is the Dolbeault cohomology of
the holomorphic vector bundle $\cV$,
$H^i(\tSi,\cO(\cV))$ is the sheaf cohomology of the sheaf $\cO(\cV)$,
and $\check{H}^i(\{ U_\alpha \}, \cO(\cV))$ is
the \v{C}ech cohomology with coefficient in the sheaf
$\cO(\cV)$ for  a good cover $\{ U_\alpha \}$ of $\tSi$.
Let $H^i(\tSi,\cV)$ denote any of the three cohomology groups in \eqref{eqn:cohomology}.

Let $\tbN_{\tau_0(\mu)}\to \cC_{\tau_0(\mu)}$ be the normal bundle
of $\cC_{\tau_0(\mu)}$ in $\cC(\tau^*E)$. Then
$\tau_\cC(\cE)= \tau_\cC(\D_1)\oplus\cdots\oplus \tau_\cC(\D_r)$,  and
$$
(\tbN_{\tau_0(\mu)})_{\tau_\cC(\cE)} =\bigoplus_{i<j} H^1\left( \tSi,\cHom(\tau_\cC(\D_j), \tau_\cC(\D_i))\right).
$$

Now, $\tau^*$ induces an anti-holomorphic map from the holomorphic
vector bundle $\cHom(\D_i, \D_j)$ to the anti-holomorphic
bundle $\tau^*\cHom(\D_i, \D_j)= \tau^*\cHom(\D_j^\vee, \D_i^\vee)$.
So $\tau$ induces an isomorphism of holomorphic vector bundles
$$
\cHom(\D_i, \D_j)\stackrel{\cong}{\longrightarrow} \overline{\tau^*\cHom(\D_j^\vee, \D_i^\vee)}
=\cHom(\tau^*\overline{\D_j^\vee}, \tau^*\overline{\D_i^\vee})=\cHom(\tau_\cC(\D_j), \tau_\cC(\D_i)).
$$

Given a local holomorphic section of $s$ of $\cHom(\D_i, \D_j)|_U$,
where $U$ is an open subset of $\tSi$, we let $\tau(s)$ denote the local holomorphic section
of $\overline{\tau^*\cHom(\D_i, \D_j)}\Bigr|_{\tau(U)}$ defined by $\tau(s)(z) = \overline{s(\tau(z))}$.
Then $\tau$ defines a conjugate linear map between
\v{C}ech complexes associated to $\cHom(\D_i,\D_j)$ and $\overline{\tau^*\cHom(\D_i,\D_j)}=\cHom(\tau_\cC(\D_j),\tau_\cC(\D_i))$,
and this in turn induces a conjugate linear map
$$
\tau: H^1(\tSi, \cHom(\D_i,\D_j))\to H^1\left(\tSi,\cHom(\tau_\cC(\D_j), \tau_\cC(\D_i))  \right).
$$
The direct sum of these maps (over $i<j$) is the conjugate linear map
$(\bN_\mu)_\cE \to (\bN_{\tau_0(\mu)})_{\tau_\cC(\cE)}$ induced by
$\tau_\cC:\cC(E)\to \cC(\tau^*E)$.

\subsection{Degree zero case}

We now specialize to the degree zero case (see \cite[Section 7]{HL1} for details).

Let $\Si^\ell_0$ be a Riemann surface of genus $\ell\geq 0$.
Let $\Si^\ell_1$ be the connected sum of $\Si^\ell_0$ and $\RP^2$, and
let $\Si^\ell_2$ be the connected sum of $\Si^\ell_0$ and a Klein bottle.
Any closed connected surface is of the form $\Si^\ell_i$, where
$\ell\geq 0$ and $i=0,1,2$. $\Si^\ell_i$ is orientable if and only if
$i=0$. Define
\begin{eqnarray*}
I_n^0 &=&\Bigl\{ \mu=(\nu,\underbrace{0,\ldots,0}_{n_0},\tau_0(\nu))\Bigl|
 \nu=\Bigl(\underbrace{\frac{k_1}{n_1},\ldots,\frac{k_1}{n_1}}_{n_1},\ldots,
\underbrace{\frac{k_r}{n_r},\ldots,\frac{k_r}{n_r} }_{n_r}\Bigr) \in I_{n',k},\\
&& \quad\quad\quad n_0>0, \quad 2n'+n_0=n, \quad \frac{k_1}{n_1}>\cdots >\frac{k_r}{n_r}>0 \Bigr\},\\
I_n^{i,\pm} &=& \Bigl\{ \mu=(\nu,\tau_0(\nu)) \Bigl|
 \nu=\Bigl(\underbrace{\frac{k_1}{n_1},\ldots,\frac{k_1}{n_1}}_{n_1},\ldots,
\underbrace{\frac{k_r}{n_r},\ldots,\frac{k_r}{n_r} }_{n_r}\Bigr) \in I_{n',k}\\
&& \quad\quad\quad2n'=n,\quad   (-1)^{n'i+k}=\pm 1,\quad
 \frac{k_1}{n_1}>\cdots >\frac{k_r}{n_r}>0\Bigr \},
\end{eqnarray*}
where $i=1,2$. Then $I_{n,0}^{\tau_0} = I_n^0 \cup I_n^{i,+}\cup I_n^{i,-}$.
To simplify notation,
we set
$$
\sigma(\mu) = (-1)^{n'i + k}
$$
for any $\mu\in I^0_n$.

We have $\Si =\Si^\ell_i$ for some $\ell\geq 0$ and $i=1,2$, and
$\tSi = \Si^{2\ell+i-1}_0$.
Let $E=P^{n,0}_\tSi \times_\rho \bC^n \cong \tSi\times \bC^n$.
The involution $\tau^\pm : E\maps E$ defines an isomorphism
$\phi^\pm: E\stackrel{\cong}{\maps} \tau^*E$. We have $\tau_\cC^\pm : \cC(E)\to \cC(E)$.
Suppose that $\mu\in I_{n,0}$, and
$\cC_\mu^{\tau_\cC^\pm}$ is nonempty. Then
$\mu \in I_n^0\cup I_n^{i,\pm}$, so $\mu$ is of the form
$$
 \mu=(\nu,\underbrace{0,\ldots,0}_{n_0},\tau_0(\nu)),
\quad \nu=\Bigl(\underbrace{\frac{k_1}{n_1},\ldots,\frac{k_1}{n_1}}_{n_1},\ldots,
\underbrace{\frac{k_r}{n_r},\ldots,\frac{k_r}{n_r} }_{n_r}\Bigr)\in I_{n',k},
$$
where $n_0\geq 0$. There exist $C^\infty$ subbundles
$D_0, \ldots, D_r$ of $E$ such that
\begin{enumerate}
\item For $i=0, \ldots, r$,
$\rank_\bC D_i = n_i,\quad \deg D_i = k_i$, where $k_0=0$.
\item
$E = D_1 \oplus \cdots D_r \oplus D_0 \oplus \tau^* D_r \oplus \cdots \oplus \tau^* D_1$.
\item $\tau^\pm$ preserves $D_0$ and switches $D_i$ with $\tau^*D_i$ for $i=1,\ldots,r$.
\end{enumerate}

Let
$$
E_\mu^0 = D_1 \oplus \cdots D_r \oplus D_0 \oplus \tau^* D_r \oplus \cdots \oplus \tau^* D_1,
$$
and define $\cB_\mu^0$ as in Section \ref{sec:levi}. Then
$\tau_\cC^\pm$ acts on $\cB_\mu^0$ by
\begin{eqnarray*}
&& \D_1 \oplus \cdots \D_r \oplus \D_0 \oplus \D_{-r}\oplus \cdots \oplus \D_{-1}\\
&\mapsto&  \tau_\cC(\D_{-1})\oplus\cdots \oplus \tau_\cC(\D_{-r})
\oplus \tau_\cC^{\pm\sigma(\mu)}(\D_0) \oplus \tau_\cC(\D_r)\oplus
\cdots \oplus \tau_\cC(\D_1),
\end{eqnarray*}
Let $\cC_{ss}(D_i)\subset \cC(D_i)$ be the semistable stratum. Any element
in the fixed locus $(\cB_\mu^0)^{\tau_\cC^\pm}$ is of the form
$$
\D_1\oplus \cdots \oplus \D_r \oplus \D_0 \oplus \tau_\cC(\D_r) \oplus \cdots \tau_\cC(\D_1)
$$
where $\D_i\in \cC_{ss}(D_i)$ for $i=1,\ldots,r$, and $\D_0\in \cC_{ss}(D_0)^{\tau_\cC^{\pm\sigma(\mu)} }$.
Therefore,
$$
(\cB_\mu^0)^{\tau_\cC^\pm}
\cong \cC_{ss}(D_0)^{\tau_\cC^{\pm\sigma(\mu)} } \times\prod_{i=1}^r \cC_{ss}(D_i).
$$

Now, $\tau_\cC^\pm$ acts on
$$
\Aut(E_\mu^0) = \Aut(D_1)\times \cdots \Aut(D_r) \times \Aut(D_0)\times
\Aut(\tau^*D_r) \times \cdots \Aut(\tau^* D_1)
$$
by
$$
(u_1,\ldots,u_r, u_0, u_{-r},\ldots, u_{-1}) \hspace{2.25in}
$$
$$
\hspace{.75in}\mapsto (\tau_\cC(u_{-1}), \ldots, \tau_\cC(u_{-r}),
\tau_\cC^{\pm\sigma(\mu)}(u_0), \tau_\cC(u_r),\ldots, \tau_\cC(u_1)
).
$$
To simplify notation, we write $\cG_\mu = \Aut(E_\mu^0)^{\tau_\cC^\pm}$.  Then we have
$$
\cG_\mu
\cong \Aut(D_0)^{\tau_\cC^{\pm\sigma(\mu)} }\times\prod_{i=1}^r \Aut(D_i).
$$

Let $\cA_\mu\subset \cA(P^{n,\pm}_\Si)$ be the equivariant Morse stratum that
corresponds to $\cC_\mu^{\tau^\pm_\cC} \subset \cC(E)^{\tau^\pm_\cC}$.
As in the orientable case, the inclusion $i: (\cB^0_\mu)^{\tau_\cC^\pm} \hookrightarrow \cC_\mu^{\tau_\cC^\pm}$
induces a weak homotopy equivalence
$$
\left( (\cB^0_\mu)^{\tau_\cC^\pm} \right)_{h\cG_\mu }
\xrightarrow{\sim} \Bigl( \cC_\mu^{\tau_\cC^\pm}\Bigr)_{h (\Aut(E)^{\tau_\cC^\pm})}.
$$
We now have a sequence of maps (where $\sim$ denotes a weak homotopy equivalence)
\begin{eqnarray*}
\left(\cA_\mu\right)_{h \Aut(P^{n,\pm}_\Si)}
\cong \bigl(\cC_\mu^{\tau_\cC^\pm} \bigr)_{h (\Aut(E,h)^{\tau_\cC^\pm})}
\xleftarrow{\sim}  \bigl(\cC_\mu^{\tau_\cC^\pm} \bigr)_{h (\Aut(E)^{\tau_\cC^\pm})} \\
\xleftarrow{\sim} \left( (\cB^0_\mu)^{\tau_\cC^\pm} \right)_{h\cG_\mu } 
\isom \left(\cC_{ss}(D_0)^{\tau_\cC^{\pm\sigma(\mu)}}\right)_{h\bigl(\Aut(D_0)^{\tau_\cC^{\pm\sigma(\mu)} }\bigr) }
\times\prod_{j=1}^r \cC_{ss}(D_i)_{h\Aut(D_i)}
\end{eqnarray*}
When $\mu\in I_n^{i,\pm}$, we do not have the first factor
$\cC_{ss}(D_0)$.

Let $\tbN_\mu$ be the normal bundle of $\cC_\mu$ in $\cC(E)$.
Given
$$
\cE=\D_1\oplus \cdots\oplus \D_r \oplus \D_0 \oplus \tau_\cC(\D_r) \oplus
\cdots \oplus \tau_\cC(\D_1) \in (\cB^0_\mu)^{\tau_\cC^\pm},
$$
we have
\begin{eqnarray*}
&& (\tbN_\mu)_\cE = H^1(\tSi,\cEnd''(\cE))\\
&=& \bigoplus_{0<i<j} H^1\left(\tSi, \cHom(\D_i,\D_j)\right) \oplus
\bigoplus_{0<i<j} H^1\left(\tSi, \cHom(\tau_\cC(\D_j), \tau_\cC(\D_i)\right) \\
&& \oplus \bigoplus_{0<i,j}H^1\left(\tSi, \cHom(\D_i, \tau_\cC(\D_j))\right)\\
&& \oplus \bigoplus_{i>0}H^1\left(\tSi,\cHom(\D_i,\D_0)\right) \oplus
\bigoplus_{i>0}H^1\left(\tSi, \cHom(\D_0, \tau_\cC(\D_i))\right).
\end{eqnarray*}

By the discussion of Section \ref{sec:levi},
$\tau$ induces conjugate linear maps of
complex vector spaces:
\begin{eqnarray*}
H^1(\tSi,\cHom(\D_i,\D_j)) & \to&  H^1(\tSi, \cHom(\tau_\cC(\D_j), \tau_\cC(\D_i))),
\textup{ and its inverse},\\
H^1(\tSi,\cHom(\D_i, \tau_\cC(\D_j)) & \to & H^1(\tSi, \cHom(\D_j, \tau_\cC(\D_i)),\\
H^1(\tSi,\cHom(\D_i, \D_0))&\to &H^1(\tSi, \cHom(\D_0,\tau_\cC(\D_i))), \textup{ and its inverse}.
\end{eqnarray*}

Let $\bN_\mu$ be the normal bundle of $\cA_\mu$ in $\cA(P^{n,\pm}_\Si)$, or equivalently,
the normal bundle of $\cC_\mu^{\tau_\cC^\pm}$ in $\cC(E)^{\tau_\cC^\pm}$. Then
\begin{equation}\label{eqn:Nmu}
\begin{aligned}
& (\bN_\mu)_\cE = H^1(\tSi,\cEnd''(\cE))^\tau\\
\cong & \bigoplus_{0<i<j} H^1\left(\tSi, \cHom(\D_i,\D_j)\right) \oplus
\bigoplus_{0<i<j} H^1\left(\tSi,\cHom(\D_i,\tau_\cC(\D_j))\right) \\
& \oplus \bigoplus_{i>0}H^1\left(\tSi, \cHom(\D_i, \tau_\cC(\D_i))\right)^\tau
\oplus \bigoplus_{i>0}H^1\left(\tSi,\cHom(\D_i,\D_0)\right)
\end{aligned}
\end{equation}

Let $i: (\cB^0_\mu)^{\tau_\cC^\pm} \hookrightarrow \cC_\mu^{\tau_\cC^\pm}$
denote the inclusion map.
By \eqref{eqn:Nmu},
$i^*\bN_\mu = \bN_\mu^\bC \oplus \bN_\mu^\bR$, where
\begin{eqnarray*}
(\bN_\mu^\bC)_\cE &=& \bigoplus_{0<i<j} H^1\left(\tSi, \cHom(\D_i,\D_j)\right) \oplus
\bigoplus_{0<i<j} H^1\left(\tSi,\cHom(\D_i,\tau_\cC(\D_j))\right) \\
&& \oplus\bigoplus_{i>0}H^1\left(\tSi,\cHom(\D_i,\D_0)\right);\\
(\bN_\mu^\bR)_\cE &=&\bigoplus_{i>0}H^1\left(\tSi, \cHom(\D_i, \tau_\cC(\D_i))\right)^\tau.
\end{eqnarray*}

Note that $(\bN_\mu^\bC)_{h\cG_\mu}\to \left( (\cB_\mu^0)^{\tau_\cC^\pm}\right)_{h\cG_\mu}$ 
is a complex vector bundle, thus
an oriented real vector bundle. Hence orientability of
$i^*\bN_\mu$ is equivalent to orientability of the real vector bundle
$$
(\bN_\mu^\bR)_{h\cG_\mu} 
\to \left( (\cB_\mu^0)^{\tau_\cC^\pm} \right)_{h\cG_\mu}
=\left(\cC_{ss}(D_0)^{\tau_\cC^{\pm\sigma(\mu)}}\right)_{h\bigl(\Aut(D_0)^{\tau_\cC^{\pm\sigma(\mu)} }\bigr) }
\times \prod_{i=1}^r \cC_{ss}(D_i)_{h\Aut(D_i)}.
$$
We have
$$
\bN_\mu^\bR=\bigoplus_{i=1}^r V_i,\quad
(V_i)_\cE = H^1(\tSi,\cHom(\D_i,\tau_\cC(\D_i)) )^\tau.
$$

Let $D$ be a rank $n$, degree $k>0$ complex
vector bundle over $\tSi$.
Let $\mathbb{V}_{n,k}$  be the
$\Aut(D)$-equivariant real vector bundle
over $\cC_{ss}(D)$ whose fiber at $\D$ is
$H^1(\tSi,\cHom(\D, \tau_\cC(\D))^\tau$.

The following result will be a direct consequence of
Lemma \ref{X-zero} (the $\tg=0$ case) and Theorem~\ref{thm:X}
(the $\tg>0$ case) in Section \ref{rep-var} below.
\begin{tm}$\label{Css}$
Let $D$ be a rank $n$, degree $k>0$
complex vector bundle over a Riemann surface
of genus $\tg$. If $n=1$ or $\tg\neq 1$ then
$\left(\mathbb{V}_{n,k}\right)_{h\Aut(D)} \to \cC_{ss}(D)_{h \Aut(D)}$
is orientable.
\end{tm}

Suppose that $\Si$ is diffeomorphic to the Klein bottle,
so that its orientable double cover $\tSi$ is a Riemann
surface of genus $\tg=1$. Note that
$$
I_2^{\tau_0} =\{ (k,-k)\mid k\in \bZ_{\geq 0}\},\quad
I_3^{\tau_0}= \{ (k,0,-k)\mid k\in \bZ_{\geq 0}\}.
$$
The $k=0$ case corresponds to open strata
whose normal bundles are of rank zero. From the above discussion,  when $\mu=(k,-k)$ or $(k,0,-k)$,
where $k>0$, we have
$$
(\bN_\mu^\bR)_\cE = H^1(\tSi,\cHom(\D_1,\tau_\cC(\D_1)))^\tau
$$
where $\D_1$ is a rank 1, degree $k$ holomorphic bundle
over $\tSi$.  Therefore Theorem \ref{Css} implies Theorem \ref{thm:main},
our main orientability theorem.

\subsection{Reduction to representation varieties}$\label{rep-var}$

We consider the following equivariant real vector bundles:
\begin{enumerate}
\item The $\Aut(E)$-equivariant vector bundle
$\mathbb{V}_{n,k} \to \cC_{ss}(E)$.
\item The $\Aut(E,h)$-equivariant vector bundle
$\mathbb{V}_{n,k}\to \cC_{ss}(E)$, or equivalently,
the $\tilde{\cG}$-equivariant vector bundle
$\mathbb{V}_{n,k} \to \cA_{ss}(P^{n,k}_\tSi)$,
where $\tilde{\cG}= \Aut(P^{n,k}_\tSi)$ and
$\cA_{ss}(P^{n,k}_\tSi)$ is the open
Morse stratum.
\item The $\tilde{\cG}$-equivariant vector bundle
$i_{n,k}^* \mathbb{V}_{n,k}\to \cN_{ss}(P^{n,k}_\tSi)$,
where $i_{n,k}:\cN_{ss}(P^{n,k}_\tSi)\hookrightarrow
\cA_{ss}(P^{n,k}_\tSi)$ is the inclusion of the space of central
Yang-Mills connections on $P^{n,k}_\tSi$.
\end{enumerate}

The inclusion $\Aut(E,h)\subset \Aut(E)$ is a homotopy
equivalence; $\cA_{ss}(P^{n,k}_\tSi)$ is
the stable manifold of $\cN_{ss}(P^{n,k}_\tSi)$,
and the gradient flow of the Yang-Mills functional
gives a $\tilde{\cG}$-equivariant
deformation retraction $\cA_{ss}(P^{n,k}_\tSi) \to \cN_{ss}(P^{n,k}_\tSi)$.
Therefore $\left(\mathbb{V}_{n,k}\right)_{ h\Aut(E)} \to \cC_{ss}(E)_{h\Aut(E)}$ is
orientable if and only if
$\left(i^*_{n,k}\mathbb{V}_{n,k}\right)_{h \tilde{\cG}}
\to \cN_{ss}(P^{n,k}_\tSi)_{h\tilde{\cG}}$ is orientable.

We fix a base point $x_0\in \tSi$, and let
$ev: \tilde{\cG}\to U(n)$ be the evaluation
at $x_0$. Then $ev$ is a surjective group
homomorphism, and the kernel $\tilde{\cG}_0$ is
the based gauge group. Therefore, $\tilde{\cG}_0$ is
a normal subgroup of $\tilde{\cG}$, and
$\tilde{\cG}/\tilde{\cG}_0 = U(n)$.
The group $\tilde{\cG}_0$ acts freely on $\cN_{ss}(P^{n,k}_\tSi)$.
Let $\tg$ be the genus of $\tSi$. The representation variety
of central Yang-Mills connections on $P^{n,k}_\tSi$ is given
by
$$
\ymU{\tg}{0}_{\kn}
=\{ V \in U(n)^{2\tg}\mid \fm(V)= e^{-2\pi\sqrt{-1}k/n} I_n \}
$$
where $\fm(a_1, b_1,\ldots, a_{\tg}, b_{\tg})=\prod_{i=1}^{\tg} [a_i,b_i]$. (See \cite[Section 6.1]{HL1} for
the definition of  $\ymU{\tg}{0}_\mu$ for a general
Atiyah-Bott type $\mu$.) There is a homeomorphism
$$
\cN_{ss}(P^{n,k}_\tSi) /\tilde{\cG}_0 \cong \ymU{\tg}{0}_{\kn}.
$$
The $\tilde{\cG}$-equivariant vector bundle
$i^*_{n,k}\mathbb{V}_{n,k}$ over $\cN_{ss}(P^{n,k}_\tSi)$ descends
to a $U(n)$-equivariant vector bundle $V_{n,k}$ over
$\ymU{\tg}{0}_{\kn}$. Therefore orientability of
$\left(i^*_{n,k}\mathbb{V}_{n,k}\right)_{h \tilde{\cG}} \to
\cN_{ss}(P^{n,k}_\tSi)_{h\tilde{\cG}}$ will follow from orientability
of the bundle
${\left(V_{n,k}\right)}_{hU(n)}\to \left(\ymU{\tg}{0}_{\kn}
\right)_{h U(n)}$.

We will need the following lemma regarding orientability
of equivariant vector bundles.

\begin{lm}$\label{equiv-or}$
Let $G$ be a compact, connected Lie group and let $X$ be a
paracompact $G$-space.  Then a $G$-equivariant real vector bundle $W\to X$ is orientable
if and only if the vector bundle $EG\cross_G W \maps EG\cross_G X$ is orientable.
\end{lm}
\begin{proof} Since $W$ is the restriction of $W_{hG}$ to a fiber of the projection $X_{hG}\to BG$, the ``if'' direction is immediate.
Now assume $W$ is orientable.  Then $\det W \to X$ is
trivial, and it will suffice to show that
$$\det (EG\cross_G W) \isom EG\cross_G (\det W) \maps EG\cross_G X$$
is trivial.  Since $X$ is paracompact and $G$ is compact, we may choose
a $G$-equivariant metric on $\det W$. The set $\det (W)_1$ of length-one vectors in $\det (W)\isom X\cross \bR$ is homeomorphic to
$X\coprod X$, so there is a section $s: X\to \det (W)$ with
image in $\det (W)_1$.

We claim that $s$ is $G$-equivariant.  Fix $x\in X$, $g\in G$.  Since $G$ is connected, there exists a
path $g_t$ from $g$ to $e$, yielding paths $s(g_t \cdot x)$ and $g_t \cdot s(x)$ from $s(g\cdot x)$ and
$g\cdot s(x)$ to $s(x)$.  By $G$-invariance of the metric, these paths lie in $\det (W)_1$, so $s(g\cdot x)$
and $g\cdot s(x)$ lie in the same path component of $\det (W)_1 \isom X\coprod X$.  Since both points are in
the fiber over $g\cdot x$, we have $s(g\cdot x) = g\cdot s(x)$.  The map $EG\cross X\to EG\cross_G \det (W)$
given by $(e, x)\mapsto [e, s(x)]$ now factors through $EG\cross_G X$, giving a nowhere-zero section of this line bundle.
\end{proof}

Since $U(n)$ is compact and connected,  to show that
$$
{\left(V_{n,k}\right)}_{h U(n)} \to
\Bigl(\ymU{\tg}{0}_{\kn}\Bigr)_{h U(n)}
$$
is orientable it suffices, by Lemma~\ref{equiv-or}, to show that
$$
V_{n,k} \to  \ymU{\tg}{0}_{\kn}
$$
is orientable.

When $\tg=0$, the definition of $\ymU{\tg}{0}_\kn$ degenerates (the reader may wish to compare
with the general definition given in \cite[Section 6.1]{HL1})
and we find that $\ymU{0}{0}_\kn$ is a single point when $k$ is a multiple of $n$, and is empty otherwise.
\begin{lm}\label{X-zero}
$$
\ymU{0}{0}_{\kn}=\begin{cases}
\{ V\in U(n)^0 | V = e^{-2\pi\sqrt{-1}k/n}I_n \} = \{I_n\}, & \frac{k}{n}\in \bZ \\
\emptyset, & \frac{k}{n}\notin \bZ
\end{cases}
$$
So $V_{n,k}\to \ymU{0}{0}_{\kn}$ is orientable whenever
$\ymU{0}{0}_\kn$ is nonempty.
\end{lm}

We will prove the following.
\begin{tm} \label{thm:X}
Let $k>0$.  The real vector bundle
$V_{n,k}\to \ymU{\tg}{0}_{\kn}$ is orientable
when $n=1$ or $\tg\geq 2$.
\end{tm}

We now explain how exactly we deduce Theorem~\ref{thm:main} from
Lemma \ref{X-zero}, Theorem~\ref{thm:X} and the previous results and arguments in this
section. We must be careful due to the fact that many of the maps we
have been considering are only weak homotopy equivalences. In
particular, if $f: X\to Y$ is a weak homotopy equivalence and $V\to
Y$ is a real vector bundle, orientability of $f^*(V)$ does not
necessarily imply orientability of $V$ (although by the Bundle
Homotopy Theorem \cite[Section 4.9]{Husemoller}, this implication does hold for homotopy
equivalences).

We need to prove orientability of the normal bundle $\left(\bN_\mu\right)_{h \cG}$ over $\left(\cA_\mu\right)_{h\cG}$.  Letting
$i:\cN_\mu \hookrightarrow \cA_\mu$ denote the inclusion of the critical set, this bundle is isomorphic
to the pullback of the bundle $i^* \left(\left(\bN_\mu\right)_{h \cG}\right)$ under the retraction $r: \cA_\mu \to \cN_\mu$ provided
by the Yang-Mills flow (because $i$ and $r$ are homotopy inverses).  Moreover, $i^* \left(\left(\bN_\mu\right)_{h \cG}\right)$ is the pull back of
a bundle $W_{hU(n)}$ over the representation variety $\left(\cN_\mu/\bG\right)_{hU(n)}$ (this reduction to
the representation variety is analogous to the argument in Section~\ref{rep-var}).
By the results in~\cite[Sections 6, 7]{HL1}, $\cA_\mu/\bG$ is an analytic set, and hence admits
a $U(n)$--equivariant triangulation~\cite{Il}.  Thus the homotopy orbit space is a CW complex, and to
prove Theorem~\ref{thm:main} we now just need to prove orientability of the bundle $W_{hU(n)}$ over the CW complex $\left(\cN_\mu/\bG\right)_{hU(n)}$.

The various weak equivalences exhibited in this section provide a weak equivalence
$$
\left((\cB^0_\mu)^{\tau_\cC^\pm}\right)_{h\Aut^0_\mu} \xrightarrow{\sim} (\cN_\mu/\bG)_{hU(n)}.
$$
By a standard CW approximation argument (for example, pull back over the singular complex of
$\left((\cB^0_\mu)^{\tau_\cC^\pm}\right)_{h\Aut^0_\mu}$), orientability of the bundle
$W_{hU(n)}\to (\cN_\mu/\bG)_{hU(n)}$ is implied by orientability of the pullback of this bundle to
$\left((\cB^0_\mu)^{\tau_\cC^\pm}\right)_{h\Aut^0_\mu}$; note that this pullback is just the restriction of
$\left(\bN_\mu\right)_{h\cG}$ to $\left((\cB^0_\mu)^{\tau_\cC^\pm}\right)_{h\Aut^0_\mu}$.
Finally, we have seen that orientability of this restricted bundle is implied by Theorem~\ref{Css},
which follows from Theorem~\ref{thm:X}.  In the subsequent sections, we will prove Theorem~\ref{thm:X}
by explicitly examining the restrictions of the bundle $V_{n,k}$ (see Theorem \ref{thm:X}) to loops generating the fundamental group
of $\ymU{\tg}{0}_{\kn}$.  Note that $\ymU{\tg}{0}_{\kn}$ is again an analytic set, hence triangulable
(here we do not need an equivariant triangulation, so the classical result of \L ojasiewicz~\cite{Lo} suffices).

\begin{rem}  For the main applications we have in mind (e.g. the Morse inequalities mentioned in the introduction), it is not strictly necessary to prove orientability of the normal bundles to the Yang-Mills strata; one simply needs Thom isomorphisms describing how each critical set contributes to the cohomology of the space $\cA_{h\cG}$.  We now explain how to deduce these isomorphisms without resorting to Illman's equivariant triangulability results, or even the non-equivariant result of \L ojasiewicz.

The partial ordering on the Yang-Mills strata defined by Atiyah and Bott \cite[Section 7]{ym} can be refined to a linear ordering in which the union of each initial segment is open (see~\cite{Ra1} for details).  Let $\cA_I$ denote the union of the strata in some initial segment $I$ in this ordering (so $\cA_I$ an open neighborhood of $\cA_{ss}$) and let $\cA_\mu$ be the next stratum.  Then, by excising the complement of a gauge-invariant tubular neighborhood (see~\cite{Ra1}) and applying the Thom Isomorphism Theorem to the (orientable) normal bundle $(\bN_\mu)_{h\cG}$ , one obtains isomorphisms
\begin{equation}\label{Thom}H^*_\cG (\cA_I\cup \cA_\mu, \cA_I) \isom H^*_\cG (\bN_\mu, \left(\bN_\mu\right)_0) \isom H^{*-c(\mu)}_\cG (\cA_\mu),
\end{equation}
where $\left(\bN_\mu\right)_0$ denotes the complement of the zero section and $c(\mu)$ is the dimension of $\bN_\mu$.  The isomorphism between the first and third terms is what we need in order to compute equivariant cohomology.

Rather than applying the Thom Isomorphism directly to $\left(\bN_\mu\right)_{h\cG}$, one may instead pull back over a CW approximation $f:X \to \left(\cA_\mu\right)_{h\cG}$.  Since $f$ is a weak equivalence and both $\bN_\mu$ and the complement of its zero section fiber over $\cA_\mu$, we have an isomorphism
$$H^*_{h\cG} (f^*\bN_\mu, \left(f^*\bN_\mu\right)_0) \isom H^*_{h\cG} (\bN_\mu, \left(\bN_\mu\right)_0).$$
To establish an isomorphism between the first and third terms in (\ref{Thom}), we need only deduce orientability of $f^*(\bN_\mu)$.  This follows from Theorem \ref{thm:X} by applying CW approximations throughout the previous argument; in fact we only need to know that the bundle $V_{n,k}$ in Theorem \ref{thm:X} is orientable after pulling back over a CW approximation $\alpha: K\stackrel{\heq}{\to} \ymU{\tg}{0}_{\kn}$.  We will show in subsequent sections that $V_{n,k}$ is orientable along loops $\{\gamma_i\}$ generating the fundamental group of $\ymU{\tg}{0}_{\kn}$.  Choosing $\gamma'_i: S^1 \to K$ such that $\alpha \circ \gamma'_i \heq \gamma_i$, the Bundle Homotopy Theorem \cite[Section 4.9]{Husemoller} implies that $\alpha^* V_{n,k}$ is orientable along the loops $\gamma'_i$, which generate $\pi_1 K$.  Since $K$ is a CW complex, this implies
(see Remark \ref{or-loops}) that $\alpha^* V_{n,k}$ is orientable, as desired.
In this approach, we do not need to use the fact that $\ymU{\tg}{0}_{\kn}$ is triangulable.
\end{rem}

\section{Fundamental Groups}\label{sec:pione}

Our orientability argument requires a calculation of fundamental
groups.

\begin{pro} \label{thm:Xdet}
For $\tg\geq 2$,
the map $\det$ induces an isomorphism
$$
\pi_1\left( \ymU{\tg}{0}_\kn \right)\stackrel{\det_*}{\longrightarrow}
\pi_1(\ymS{\tg}{0}_k)\cong \bZ^{2\tg}.
$$
\end{pro}
\begin{proof} We may assume $n\geqs 2$.
We first introduce some notation.
Let $P^{m,k}= P^{m,k}_{\tSi}$;
note that $\det(P^{n,k})=P^{1,k}$. Let $\cA(m,k) =\cA(P^{m,k})$
be the space of $U(m)$-connections on $P^{m,k}$, let
$\cG(m,k) = \Aut(P^{m,k})$ be the gauge group, and
let $\cG_0(m,k)\subset \cG(m,k)$ be the base gauge group. Let
$\cC(m,k)$ be the space of holomorphic structures
on $E^{m,k}$, the rank $m$, degree $k$ complex
vector bundle over $\tSi$. Let $\Css(m,k)\subset \cC(m,k)$
be the semi-stable stratum.

Recall that Trace$:\fu(n) \to \fu(1)$ is the derivative of
the determinant map det$: U(n) \to U(1)$ at the identity. Clearly it is
$\mathrm{ad}$-invariant and it induces a map $ \ad(P^{n,k})\to
\ad(P^{1,k})$, and thus a map $\tr:\cA(n,k) \to \cA(1,k)$.
The map $\tr$ sends a Yang-Mills $U(n)$-connection to a Yang-Mills
$U(1)$-connection. Since all Yang-Mills $U(1)$-connection are
central, the map $\tr$ descends to a map
$$\det: \ymU{\tg}{0}_\kn \to
\ymS{\tg}{0}_k$$
(recall that $\det(\exp M)=\exp(\tr M), \forall M
\in \fu(n))$. In other words, we have a commuting diagram:
$$\begin{CD}
\cN_\kn @>{\tr}>> \cN_k \\
@VV{ \hol} V @VV{\hol}V \\
\ymU{\tg}{0}_\kn @>{\det}>> \ymS{\tg}{0}_k \\
\end{CD}
$$

The determinant map $U(n)\to U(1)$ also induces a homomorphism
$\phi: \cG_0(n,k) \to \cG_0(1,k)$, and the map $\tr:
\cA(n,k) \to \cA(1,k)$ is $\phi$-equivariant.
In particular, the map $\tr: \cN_\kn \to \cN_k$ is
$\phi$-equivariant, and
 we have a well-defined map
\begin{equation*}\tr:  E\cG_0(n,k)
\cross_{\cG_0(n,k)} \cN_\kn  \maps  E\cG_0(1,k) \cross_{\cG_0(1,k)}
\cN_k
\end{equation*}
which we may identify up to homotopy with
the determinant map
$$\det:\ymU{\tg}{0}_\kn \to \ymS{\tg}{0}_k.$$ Moreover, the
Yang-Mills flow provides a gauge-equivariant deformation retraction
from the space $\Css(m,k)$ of semi-stable bundles to the critical
set $\cN_\km$~\cite{Rade}, so it suffices to show that the map
\begin{equation}\label{orbits}
\pi_1 \left( E\cG_0(n,k) \cross_{\cG_0(n,k)} \Css(n,k) \right)
\stackrel{\tr_*}{\maps} \pi_1\left( E\cG_0(1,k) \cross_{\cG_0(1,k)} \Css(1,k) \right)
\end{equation}
is an isomorphism.

We have an induced map of fibration sequences
\begin{equation}\label{fib}
\xymatrix{
    \Css(n,k) \ar[r] \ar[d]
    & E\cG_0(n,k) \cross_{\cG_0(n,k)} \Css(n,k) \ar[r] \ar[d]
    & B\cG_0(n,k) \ar[d]\\
    \Css(1,k) \ar[r] &
    E\cG_0(1,k) \cross_{\cG_0(1,k)} \Css(1,k) \ar[r]
     & B\cG_0(1,k),
}
\end{equation}
and we claim that both fibers are simply connected.
 For $n=1$, all critical connections are minimal, i.e.
there is only one stratum and thus the set of minimal Yang-Mills connections is a
deformation retraction of the total space $\cA(1,k)$, which is an affine space.
Thus, $\Css(1,k)=\cA(1,k)$ is contractible.
Since both $n$ and $\tg$ are at least 2,
the complement of $\Css(n,k)$ in the contractible space $\cA(n,k)$ may be stratified
by submanifolds of (finite) real codimension at least $2(\tg-1)(n-1)+2 \geqs 4$. Transversality arguments (as
in~\cite[Section 4]{Ra2} or~\cite{DU}) now
apply to prove simple connectivity.

Since both $\pi_1$ and $\pi_0$ of $\Css(1,k)$ and
$\Css(n,k)$ are trivial, we may now identify the map (\ref{orbits}) with
the map $\pi_1 (B\cG_0(n,k)) \stackrel{c}{\maps} \pi_1 (B \cG_0 (1, k))$
induced by diagram (\ref{fib}).
By~\cite[Section 2]{ym}, we have homotopy equivalences
$B\cG_0 (m,k) \heq \Map^{P^{m,k}}_*(\tSi, BU(m))$ for any $m$,
where  $\Map^{P^{m,k}}_*$ denotes the subspace of based maps which induce the
bundle $P^{m,k}$.
Hence we may identify the map $c$ with the determinant map
\begin{equation}\label{AB-model}
\pi_1 (\Map^{P^{n,k}}_*(\tSi, BU(n)) )\maps \pi_1 (\Map^{P^{1,k}}_*(\tSi, BU(1))).
\end{equation}
The splitting
$U(1) \to U(n)$ of $\det: U(n)\to U(1)$ induces a splitting
\begin{equation}\label{AB-split}
\Map_*(\tSi, BU(1)) \stackrel{i}{\maps} \Map_* (\tSi, BU(n))
\end{equation}
of the determinant map
$\Map_*(\tSi, BU(n))\maps \Map_*(\tSi, BU(1))$,
and hence after restricting to components (recall that $\det(P^{n,k})=P^{1,k}$)
we obtain splittings of the maps (\ref{AB-model}).  This implies that the maps (\ref{AB-model}) are surjective.

To prove that the maps (\ref{AB-model}) are also injective, it suffices to show that their domain and range are isomorphic to $\bZ^{2\tg}$.
Note that $\tSi$ is  the mapping cone of the attaching map $\eta$ for its $2$--cell, so we have a homotopy cofiber sequence $S^1 \stackrel{\eta}{\to} \bigvee_{2\tg} S^1 \to \tSi$.
For any $m\geqs 1$, applying $\Map_* (-, BU(m))$ to this sequence
gives the fibration sequence
\begin{equation}\label{map-fib}\Map_* (\tSi, BU(m)) \stackrel{r}{\maps} \Map_* (\bigvee_{2\tg} S^1, BU(m))  \stackrel{s}{\maps} \Map_* (S^1, BU(m)).\end{equation}
We have $\Map_* (\bigvee_{2\tg} S^1, BU(m)) = (\Omega BU(m))^{2\tg} \heq U(m)^{2\tg}$ and similarly $\Map_* (S^1, BU(m)) \heq U(m)$, so the fundamental groups of these spaces are $\bZ^{2\tg}$ and $\bZ$, respectively.  Since the attaching map $\eta$ can be written as a product of commutators, so can the induced map
$$s_*: \pi_1 \, \Map_* (\bigvee_{2\tg} S^1, BU(m))  \maps \pi_1\, \Map_* (S^1, BU(m)).$$
Since these groups are abelian, we see that $s_* = 0$.

Now, a classifying map for $P^{m,k}$ gives each space in (\ref{map-fib}) a basepoint, and the resulting long exact sequence in homotopy is, in part,
$$\pi_2 \Omega BU(m) = 0 \maps \pi_1 \, \Map^{P^{m,k}}_* (\tSi, BU(m)) \stackrel{r_*}{\maps} \bZ^{2\tg} \xrightarrow{s_* = 0} \bZ.$$
Hence $r_*$ is an isomorphism, which completes the proof.
\end{proof}

\section{Symmetric Representation Varieties} \label{sec:symm-rep}

In Section \ref{sec:reduction}, we reduced our main theorem (Theorem
\ref{thm:main}) to the orientability of a real vector bundle
$V_{n,k}$ over the representation variety
$$\ymU{2\ell+i-2}{0}_\kn$$
of the central Yang-Mills $U(n)$-connection on the orientable double
cover $\Si^{2\ell+i-1}_0$ of the nonorientable surface $\Si^\ell_i$
(Theorem \ref{thm:X}). In this section, we will use Proposition
\ref{thm:Xdet} to write down:
\begin{itemize}
\item[(i)] loops in $\ymU{2\ell+i-2}{0}_\kn$ that
generate the fundamental group of $\ymU{2\ell+i-1}{0}_\kn$, and
\item[(ii)] lifts of these loops under the  surjective continuous map
$$\Phi^{\ell,i}: \ZymU{i}_\kn\maps \ymU{2\ell+i-1}{0}_\kn$$
from the symmetric representation variety.

\end{itemize}

To prove  Theorem \ref{thm:X}, it suffices to examine the
orientability of the restrictions of the pull back bundle $W_{n,k}=
(\Phi^{\ell,i})^* V_{n,k}$ to the loops in (ii). This will be
carried out in Section \ref{sec:alongloops}.

\subsection{Review of symmetric representation varieties}
We recall definitions and some properties of symmetric
representation varieties introduced in \cite{HL1}.

Given $V=(\ab)\in U(n)^{2\ell}$, let $\fm(V)=\pab$. For integers
$k,n$, where $n>0$, we introduce symmetric representation varieties:
\begin{eqnarray*}
\ZymU{1}_\kn &=&
\bigl\{(V,c,V',c',-2\sqrt{-1}\pi \frac{k}{n} I_n)\mid V,V' \in U(n)^{2\ell},\ c,c'\in U(n),\\
&& \quad \fm(V) =e^{-\pi\sqrt{-1}k/n}I_n c c',\
\fm(V')=e^{\pi \sqrt{-1}k/n}I_n c' c \bigr\}\\
\ZymU{2}_\kn&=&
\bigl\{(V,d,c,V',d',c',-2\sqrt{-1}\pi \frac{k}{n} I_n)\mid V,V'\in U(n)^{2\ell},\\
&&\quad d,c,d', c'\in U(n), \fm(V) =e^{-\pi\sqrt{-1}k/n}I_n cd' c^{-1}d,\ \\
&&\quad \fm(V')=e^{\pi\sqrt{-1}k/n}I_n c' d (c')^{-1}d' \bigr\}
\end{eqnarray*}
In particular, we have homeomorphisms
\begin{eqnarray*}
\ZymS{1}_k &=& \bigl\{(V,c,V',c',-2\sqrt{-1}\pi k)\mid V,V\in U(1)^{2\ell},\\
&&\quad  c,c'\in U(1), cc'=(-1)^k \bigr\}\ \cong\ U(1)^{4\ell+1}\\
\ZymS{2}_k &=& \bigl\{(V,d,c,V',d',c',-2\sqrt{-1}\pi k)\mid V,V'\in U(1)^{2\ell},\\
&&\quad  d,c,d',c'\in U(1),\ dd'=(-1)^k \bigr\} \ \cong \
U(1)^{4\ell+3}
\end{eqnarray*}

Given $g\in U(n)$ and $V=(\ab)\in U(n)^{2\ell}$, let
$$
g V g^{-1}=(g a_1 g^{-1}, g b_1 g^{-1},\ldots, g a_\ell g^{-1}, g
b_\ell g^{-1}).
$$
With this notation, $U(n)^2$ acts on $\ZymU{i}_\kn$ by
\begin{eqnarray*}
&& (g_1,g_2)\cdot (V,c,V',c',,-2\sqrt{-1}\pi \frac{k}{n}I_n) \\
&=&  (g_1 V g_1^{-1}, g_1 c g_2^{-1},g_2 V' g_2^{-1}, g_2 c' g_1^{-1}, -2\sqrt{-1}\pi\frac{k}{n} I_n)\\
&& (g_1,g_2)\cdot (V,d,c,V',d',c',-2\sqrt{-1}\pi \frac{k}{n}I_n)\\
&=& (g_1 V g_1^{-1},g_1 d g_1^{-1}, g_1 c g_2^{-1}, g_2
V'g_2^{-1},g_2d'g_2^{-1}, g_2 c' g_1^{-1}, -2\sqrt{-1}\pi
\frac{k}{n}I_n)
\end{eqnarray*}
Define $\Phi^{\ell,i}: U(n)^{2(2\ell+i)}\times \fu(n) \to
U(n)^{2(2\ell+i-1)}\times \fu(n)$ by
\begin{eqnarray*}
\Phi^{\ell,1}(V,c,V',c',X) &=&(V, c\fr(V') c^{-1}, X)\\
\Phi^{\ell,2}(V,d,c,V',d',c',X)&=& (V, d^{-1}c\fr(V')c^{-1}d,
d^{-1}, cc',X)
\end{eqnarray*}
where $\fr(\ab)=(b_\ell,a_\ell,\ldots,b_1,a_1)$. Then
$$
\Phi^{\ell,i}(\ZymU{i}_\kn)=\ymU{2\ell+i-1}{0}_\kn.
$$

\subsection{Maps and vector bundles}

In this subsection, $i=1,2$, and $n,k$ are positive integers.

Given a rank $n$, degree $k$ holomorphic vector bundle $\D$ over
$\tSi$, $\tau_\cC(\D)$ is a rank $n$, degree $-k$ holomorphic vector
bundle over $\tSi$, and $\cHom(\D,\tau_\cC(\D)) = \D^\vee\otimes
\tau_\cC(\D)$ is a degree $-2k$, rank $n^2$ holomorphic vector
bundle over $\tSi$. The map $\D\mapsto \tau_\cC(\D)$ defines
\begin{equation}\label{eqn:Ztau}
\tau: \ZymU{i}_\kn \to \ZymU{i}_{\mkn}.
\end{equation}
The map $\D\mapsto \cHom(\D,\tau_\cC(\D))$ defines
\begin{equation}\label{eqn:Zphi}
\phi:\ZymU{i}_\kn \to Z^{\ell,i}_{\mathrm{YM}}(U(n^2))_\mtwokn.
\end{equation}
The map $\cM \mapsto \tau_\cC(\cM^\vee)=\overline{\tau^*\cM}$
defines
\begin{equation}\label{eqn:htau}
\hat{\tau}: Z^{\ell,i}_{\mathrm{YM}}(U(n^2))_\mtwokn \to
Z^{\ell,i}_{\mathrm{YM}}(U(n^2))_\mtwokn.
\end{equation}

There is a map $U(n)\times U(n) \to U(n^2)$ given by $(A,B)\mapsto
A\otimes B$. More explicitly,
$$
(A\otimes B)_{ij, pq} = A_{ip}B_{jq},\quad 1\leq i,j,p,q \leq n.
$$
Note that $I_n\otimes I_n=I_{n^2}$. In particular, when $n=1$, this
map is the multiplication: $U(1)\times U(1)\to U(1)$,
$(c_1,c_2)\mapsto c_1 c_2$.

We introduce some notation.
\begin{enumerate}
\item[(i)] Given $A=(A_{ij})\in U(n)$, let $\bar{A}=(\bar{A}_{ij})$ be the complex
conjugate of $A$. Then $\bar{A}= (A^t)^{-1}$.

\item[(ii)] We define a complex linear involution $T$ on $\bC^n\otimes \bC^n \cong \bC^{n^2}$
by $T(u\otimes v) = v\otimes u$ for $u,v\in \bC^n$. Then
$T\in O(n^2)\subset U(n^2)$. We have
$$
T_{ij,pq}= \delta_{iq}\delta_{pj},\quad
T=T^t = T^{-1},\quad
(TCT^{-1})_{ij,pq}= C_{ji,qp}.
$$

\item[(iii)]
Define an involution
$$f_T:Z_{\mathrm{YM}}^{\ell,i}(U(n^2))_\mtwokn\to
Z_{\mathrm{YM}}^{\ell,i}(U(n^2))_\mtwokn$$ by
\begin{eqnarray*}
&& (V, c, V',c', 4\sqrt{-1}\pi\frac{k}{n} I_{n^2})
\mapsto (I_{n^2},T)\cdot (V,c, V',c', 4\sqrt{-1}\pi\frac{k}{n} I_{n^2}),\\
&&
 (V,d,c, V',d',c', 4\sqrt{-1}\pi\frac{k}{n} I_{n^2})
\mapsto (I_{n^2},T)\cdot (V,d,c, V',d',c', 4\sqrt{-1}\pi\frac{k}{n}
I_{n^2}).
\end{eqnarray*}
More explicitly,
\begin{eqnarray*}
f_T(V,c,V',c', 4\sqrt{-1}\pi\frac{k}{n}I_{n^2}) &=& (V,c T^{-1}, T
V' T^{-1}, Tc', 4\sqrt{-1}\pi\frac{k}{n}I_{n^2}),
\\f_T(V,d, c,V',d', c', 4\sqrt{-1}\pi\frac{k}{n}I_{n^2})
&=& (V,d, c T^{-1}, T V' T^{-1}, Td' T^{-1}, Tc',
4\sqrt{-1}\pi\frac{k}{n}I_{n^2}).
\end{eqnarray*}

\item[(iv)] Given $V=(\ab)$ and
$V'=(a_1',b_1',\ldots,a_\ell',b_\ell')$ in $U(1)^{2\ell}$, define
$$
V V'=(a_1 a_1',b_1 b_1',\ldots, a_\ell a_\ell', b_\ell b_\ell') \in
U(1)^{2\ell}.
$$
\item[(v)] Given $V=(\ab)$ and
$V'=(a'_1,b'_1,\ldots,a'_\ell,b'_\ell)$ in  $U(n)^{2\ell}$, define
$$
V\otimes V' = (a_1\otimes a_1', b_1\otimes b_1',\ldots,
a_\ell\otimes a'_\ell, b_\ell\otimes b_\ell')\in U(n^2)^{2\ell}.
$$
\item[(vi)] Given $V=(\ab)\in U(n)^{2\ell}$,
define
$$
\bV=(\bab) \in U(n)^{2\ell}.
$$
\end{enumerate}
It is straightforward to check that for $A,B\in U(n)$,
\begin{equation}
\overline{A\otimes B} = \bar{A}\otimes \bar{B},
\end{equation}
\begin{equation}\label{eqn:TABT}
T(A\otimes B)T^{-1} = B\otimes A.
\end{equation}
If $A, B$ are diagonal with respect to the standard basis $\{
e_i\mid i=1,\ldots, n\}$ of $\bC^n$, then $A\otimes B$ is diagonal
with respect to the basis $\{ e_i\otimes e_j \mid i,j=1,\ldots,n\}$
of $\bC^{n^2}$.

With the above notation, we  have the following explicit description
of the maps $\tau$, $\phi$, $\hat{\tau}$ in \eqref{eqn:Ztau},
\eqref{eqn:Zphi}, \eqref{eqn:htau}, respectively.

The involution $\tau: \ZymU{i}_\kn \to  \ZymU{i}_\mkn$ is given by
\begin{eqnarray*}
&& (V,c,V',c',-2\sqrt{-1}\pi \frac{k}{n} I_n) \mapsto
 (V',c', V,c,2\sqrt{-1}\pi \frac{k}{n} I_n),\quad i=1,\\
&&(V,d,c,V',d',c',-2\sqrt{-1}\pi \frac{k}{n} I_n) \mapsto (V',d',
c',V,d,c,2\sqrt{-1}\pi \frac{k}{n} I_n),\quad i=2.
\end{eqnarray*}

The map $\phi:\ZymU{i}_\kn \to
Z^{\ell,i}_{\mathrm{YM}}(U(n^2))_\mtwokn$ is given by
\begin{eqnarray*}
&& (V,c,V',c',-2\sqrt{-1}\pi \frac{k}{n} I_n) \\
&\mapsto& (\bV\otimes V', \bc\otimes c',
\bV'\otimes V, \bc'\otimes c,4\sqrt{-1}\pi \frac{k}{n} I_{n^2}),\quad i=1,\\
&& (V,d,c,V',d',c',-2\sqrt{-1}\pi \frac{k}{n} I_n)\\
&\mapsto& (\bV\otimes V', \bd \otimes d', \bc \otimes c', \bV'
\otimes V, \bd' \otimes d, \bc' \otimes c, 4\sqrt{-1}\pi \frac{k}{n}
I_{n^2}),\quad i=2.
\end{eqnarray*}
Letting $\phi_T = f_T\circ \phi$, we see that $\phi$ and $\phi_T$ define the
same map to the quotient of
$Z^{\ell,i}_{\mathrm{YM}}(U(n^2))_\mtwokn$ by $U(n^2)^2$. We have
\begin{equation}\label{eqn:RPphiT}
\begin{aligned}
& \phi_T(V,c,V',c',-2\sqrt{-1}\pi \frac{k}{n} I_n) \\
=& (\bV\otimes V', \left(\bc\otimes c'\right) T, V\otimes \bV',
\left(c\otimes \bc'\right) T ,4\sqrt{-1}\pi \frac{k}{n}
I_{n^2}),\quad i=1,
\end{aligned}
\end{equation}
\begin{equation}\label{eqn:KphiT}
\begin{aligned}
& \phi_T(V,d,c,V',d',c',-2\sqrt{-1}\pi \frac{k}{n} I_n)\\
=& ( \bV\otimes V', \bd  \otimes d', \left(\bc \otimes c'\right) T,
V\otimes \bV'  , d\otimes \bd' , \left(c\otimes \bc'\right) T,
4\sqrt{-1}\pi \frac{k}{n} I_{n^2}),\quad i=2.
\end{aligned}
\end{equation}
The involution $\hat{\tau}:Z^{\ell,i}_{\mathrm{YM}}(U(n^2))_\mtwokn
\to Z^{\ell,i}_{\mathrm{YM}}(U(n^2))_\mtwokn$ is given by
\begin{eqnarray*}
&& (V,c,V',c',4\sqrt{-1}\pi \frac{k}{n} I_{n^2}) \mapsto ( \bV' ,
\bc', \bV, \bc,
4\sqrt{-1}\pi \frac{k}{n} I_{n^2}),\quad i=1,\\
&& (V,d, c,V',d', c',4\sqrt{-1}\pi \frac{k}{n} I_{n^2}) \mapsto (
\bV', \bd',\bc', \bV, \bd, \bc, 4\sqrt{-1}\pi \frac{k}{n}
I_{n^2}),\quad i=2.
\end{eqnarray*}
Let $Z^{\ell,i}_{\mathrm{YM}}(U(n^2))_\mtwokn^{\hat{\tau}}$ be the
fixed locus of $\hat{\tau}$. Then
$$
\phi_T\Bigl(\ZymU{i}_{\kn}\Bigr) \subset
Z^{\ell,i}_{\mathrm{YM}}(U(n^2))_\mtwokn^{\hat{\tau}}.
$$

In particular, when $n=1$, we have $\bc=c^{-1}$ and $\bd=d^{-1}$.
The involution $\tau:\ZymS{i}_k \to  \ZymS{i}_{-k}$ is given by
\begin{eqnarray*}
&& (V,c,V',(-1)^k \bc,-2\sqrt{-1}\pi k) \mapsto
(V',(-1)^k \bc,V,c,2\sqrt{-1}\pi k),\quad i=1,\\
&& (V,d,c,V',(-1)^k \bd,c',-2\sqrt{-1}\pi k) \mapsto (V',(-1)^k
\bd,c',V,d,c,2\sqrt{-1}\pi k),\quad i=2.
\end{eqnarray*}
The map $\phi=\phi_T:\ZymS{i}_k \to \ZymS{i}_{-2k}$ is given by
\begin{eqnarray*}
&& (V,c,V',(-1)^k c,-2\sqrt{-1}\pi k)\\
&\mapsto& ( \bV V',(-1)^k
\bc^2,\bV' V,(-1)^k c^2,4\sqrt{-1} \pi k),\quad i=1,\\
&& (V,d,c, V',(-1)^k \bd,c',-2\sqrt{-1}\pi k) \\
&\mapsto& (\bV V',(-1)^k \bd^2, \bc c', \bV' V,(-1)^k d^2, \bc'
c,4\sqrt{-1} \pi k),\quad i=2.
\end{eqnarray*}
The involution $\hat{\tau}:\ZymS{1}_{-2k} \to \ZymS{1}_{-2k}$ is
given by
\begin{eqnarray*}
&& (V,c,V',\bc, 4\sqrt{-1}\pi k)
\mapsto ( \bV', c, \bV , \bc , 4\sqrt{-1}\pi k),\quad i=1,\\
&&(V,d, c,V', \bd , c', 4\sqrt{-1}\pi k) \mapsto ( \bV' ,d, \bc' ,
\bV , \bd , \bc, 4\sqrt{-1}\pi k),\quad i=2.
\end{eqnarray*}

\begin{df}\label{Uvec}
When $k>0$,  let $U_{n,k}\to Z_{\mathrm{YM}}^{\ell,i}(U(n^2))_{\mtwokn}$ be the
complex vector bundle whose fiber at $\cM$ is $H^1(\tSi, \cM)$,
where $\cM$ is a polystable holomorphic vector bundle of rank $n^2$,
degree $-2nk$.
\end{df}
 From the discussion in Section \ref{sec:levi}, the
involution $\tau:\tSi\to \tSi$  induces a conjugate linear map
$$
\hat{\tau}: (U_{n,k})_\cM = H^1(\tSi,\cM) \to
(U_{n,k})_{\hat{\tau}(\cM)} = H^1(\tSi,\overline{\tau^*\cM}).
$$
So the involution $\hat{\tau}$ on
$Z_{\mathrm{YM}}^{\ell,i}(U(n^2))_{\mtwokn}$ lifts to an involution
$\hat{\tau}$ on $U_{n,k}$, and the fixed locus
$U_{n,k}^{\hat{\tau}}$ is a real vector bundle over
$Z_{\mathrm{YM}}^{\ell,i}(U(n^2))^{\hat{\tau}}_{\mtwokn}$. We have
$$
\phi_T^* U^{\hat{\tau}}_{n,k} = W_{n,k},\quad \rank_\bR W_{n,k} =
\rank_\bR U_{n,k}^{\hat{\tau}} =\rank_\bC U_{n,k}= 2nk +
n^2(2\ell+i-2).
$$

\begin{df}\label{def:Uvir}
When $k=0$,  let $U_n^\vir \to Z_{\mathrm{YM}}^{\ell,i}(U(n^2))_{0,\ldots,0}
=Z_{\mathrm{flat}}^{\ell,i}(U(n^2))$ be the
virtual complex vector bundle whose fiber at $\cM$ is
$H^1(\tSi, \cM)-H^0(\tSi,\cM)$,
where $\cM$ is a polystable holomorphic vector bundle of rank $n^2$,
degree $0$.
\end{df}
The involution $\tau:\tSi\to \tSi$  induces a conjugate linear map
$$
\hat{\tau}: (U_n^\vir)_\cM =H^1(\tSi, \cM)-H^0(\tSi,\cM) \to
(U_n^\vir)_{\hat{\tau}(\cM)}=H^1(\tSi,
\hat{\tau}(\cM))-H^0(\tSi,\hat{\tau}(\cM)).
$$
$Z_{\mathrm{flat}}^{\ell,i}(U(n^2))$ lifts to an involution
$\hat{\tau}$ on $U_n^\vir$, and the fixed locus
$(U_n^\vir)^{\hat{\tau}}$ is a virtual real vector bundle over
$Z_{\mathrm{flat}}^{\ell,i}(U(n^2))^{\hat{\tau}}$ of rank
$n^2(2\ell+i-2)$.

\subsection{Loops in the symmetric representation variety}
Let
$$
\Psi: \ymU{2\ell+i-1}{0}_\kn \lra \ymS{2\ell+i-1}{0}_k
$$
be as in Proposition \ref{thm:Xdet}. In particular, $\Psi$ is the
identity map when $n=1$.

Given $V=(\ab)\in U(n)^{2\ell}$, define
$$
\det(V)=(\det (a_1), \det (b_1),\ldots, \det (a_\ell), \det (b_\ell)
) \in U(1)^{2\ell}.
$$
Then $\Psi\circ \Phi^{\ell,i}:\ZymU{i}_\kn \lra \ymS{2\ell+i-1}{0}_k
\cong U(1)^{2(2\ell+i-1)}$ is given by
\begin{equation}\label{eqn:PhiRP}
\Psi\circ\Phi^{\ell,1}(V,c,V',c',X) = \bigl(\det(V), \det(\fr(V'))
,\Tr(X)\bigr)
\end{equation}
\begin{equation}\label{eqn:PhiK}
\Psi\circ\Phi^{\ell,2}(V,d,c,V',d',c',X) = \bigl(\det(V),
\det(\fr(V')), \det(d)^{-1},\det(cc'),\Tr(X)\bigr)
\end{equation}
In the rest of this subsection, we write $\Phi$ instead of
$\Phi^{\ell,i}$.

\begin{rem} \label{or-loops}
The following observation will be useful. Let $M$ be a CW complex,
and let $E\to M$ be a real vector bundle. Then $E$ is orientable if
and only if $w_1(E)=0$, or equivalently, $\gamma \cap w_1(E)=0$ for
all $\gamma\in H_1(M;\bZ/2\bZ)$. (Recall here that orientability is equivalent, over any base,
to triviality of the determinant line bundle, and line bundles over a CW complex X are classified
by $w_1 \in H^1(X; \bZ/2\bZ) \isom [X, \RP^\infty = K(\bZ/2\bZ, 1)]$.)
Now suppose that $\pi_1(M)$ is a free
abelian group generated by loops $\gamma_1,\ldots, \gamma_r$. Then the $\gamma_i$ represent generators for $H_1 (X; \bZ/2\bZ)$, so
$E$ is orientable if and only if $[\gamma_i] \cap w_1(E)=0$ for
$i=1,...,r$, or equivalently,  $E \bigr|_{\gamma_i}$ is orientable
for $i=1,...,r$.
\end{rem}
\subsubsection{$n=1$, $i=1$} \label{sec:RPloop}
When $\ell=0$, $X^{0,0}(U(1))_{k,\ldots,k}$ consists of a single
point. We may assume $\ell\geq 1$. For $i=1,\ldots,\ell$, let
$\alpha_i$, $\beta_i$, $\alpha'_i$, and $\beta'_i$ be loops on
$\tSi$ which generate the fundamental group of $\tSi=\Si^{2\ell}_0$,
i.e., the holonomies along $\alpha_i$, $\beta_i$, $\alpha'_i$, and
$\beta'_i$ are $a_i$, $b_i$, $a'_i$, and $b'_i$ respectively. Let
$\tilde{\alpha}_i:  S^1 \to \ZymS{1}_k$ be the loop defined by
\begin{eqnarray*}
&& a_i= e^{\sqrt{-1}\theta},\quad a_j = 1, j\neq i; \quad b_j=a_j'=b_j'=1, j=1,\ldots,\ell;\\
&& c=c'=\sqrt{-1}^k, ~X=-2\sqrt{-1}\pi k.
\end{eqnarray*}
We define $\tilde{\alpha}'_i$, $\tilde{\beta}_j$, $\tilde{\beta}'_j$
similarly. Then $\Phi\circ \tilde{\alpha}_i: S^1 \to
\ymS{2\ell}{0}_k$ is a loop defined by
$$
a_i= e^{\sqrt{-1}\theta},\quad a_j = 1, j\neq i,\quad b_j=1,
j=1,\ldots,2\ell, \quad X=-2\sqrt{-1}\pi k.
$$
Thus the $4\ell$ loops
\begin{equation} \label{eqn:oneRPloop}
\Phi \circ \tilde{\alpha}_i,\ \Phi\circ\tilde{\beta}_i,\ \Phi\circ
\tilde{\alpha}'_i,\ \Phi \circ \tilde{\beta}'_i,\quad
i=1,\ldots,\ell.
\end{equation}
generate the fundamental group of $\ymS{2\ell}{0}_k\cong
U(1)^{4\ell}$.

To show that $V_{1,k} \to \ymS{2\ell}{0}_k$ is orientable, it
suffices to show that its restriction to each of the $4\ell$ loops
in \eqref{eqn:oneRPloop} is orientable, or equivalently:
\begin{pro}\label{thm:oneRP}
Let $\ell\geq 1$. The restriction of $W_{1,k}\to \ZymS{1}_k$ to each
of the following $4\ell$ loops is orientable:
$$
\tilde{\alpha}_i,\ \tilde{\beta}_i,\ \tilde{\alpha}'_i,\
\tilde{\beta}'_i,\quad i=1,\ldots,\ell.
$$
\end{pro}
\begin{proof}
See Section \ref{sec:RPproof}.
\end{proof}

\subsubsection{$n=1$, $i=2$}\label{sec:Kloop}
In this subsection, we assume that $\ell\geq 0$. For
$i=1,\ldots,\ell$, the holonomies along $\alpha_i$, $\beta_i$,
$\alpha'_i$, and $\beta'_i$ are $a_i$, $b_i$, $a'_i$, and $b'_i$
respectively. Let $\gamma$ be the curve from $p_+$ to $p_-$,
$\gamma'$ be the curve from $p_-$ to $p_+$, $\delta$ be the loop
starting at $p_+$, and $\delta'$ be the loop starting at $p_-$ i.e.,
the holonomies along $\gamma,\gamma',\delta,\delta'$ are $c, c', d,
d'$, respectively. The fundamental group of $\tSi=\Si^{2\ell+1}_0$
is generated by the $4\ell+2$ loops
$$
\alpha_i, \beta_i, \alpha'_i, \beta'_i,\quad i=1,\ldots,\ell,\quad
\delta, \gamma\gamma'.
$$
In particular, when $\ell=0$, the fundamental group of $\Si^1_0\cong
S^1\times S^1$ is generated by $\delta$ and $\gamma\gamma'$.

Let $\tilde{\alpha}_i:S^1\to \ZymS{2}_k$ be the loop defined by
\begin{eqnarray*}
&& a_i=e^{\sqrt{-1}\theta},\quad a_j=1 \textup{ if } j\neq i,\quad
b_j=a_j'=b_j'=1, j=1,\ldots,\ell,\\
&& c=c'=1,\quad d=\sqrt{-1}^{k+1}, d'=\sqrt{-1}^{k-1},\quad
X=-2\sqrt{-1}\pi k.
\end{eqnarray*}
We define $\tilde{\beta}_i$, $\tilde{\alpha}'_i$, $\tilde{\beta}'_i$
similarly.

Let $\tilde{\delta}: S^1 \to\ZymS{2}_k$ be the loop defined by
\begin{eqnarray*}
&& d=\sqrt{-1}^{k+1}e^{-\sqrt{-1}\theta},\quad
d'=\sqrt{-1}^{k-1}e^{\sqrt{-1}\theta},\quad c=c'=1\\
&& a_j=b_j=a_j'=b_j'=1,\ j=1,\ldots,\ell,\quad X=-2\sqrt{-1}\pi k.
\end{eqnarray*}

Let $\tilde{\gamma}: S^1 \to \ZymS{2}_k$ be the loop defined by
\begin{eqnarray*}
&& c=e^{\sqrt{-1}\theta},\quad c'=1,\quad d=\sqrt{-1}^{k+1},
~d'=\sqrt{-1}^{k-1},\\
&& a_j=b_j=a_j'=b_j'=1, j=1,\ldots,\ell,\quad X=-2\sqrt{-1}\pi k.
\end{eqnarray*}

Then $\Phi\circ \tilde{\alpha}_i: S^1 \to \ymS{2\ell+1}{0}_k$ is a
loop defined by
\begin{eqnarray*}
&& a_i = e^{\sqrt{-1}\theta},\quad a_j = 1\textup{ if } j\notin\{i,2\ell+1\},\quad
a_{2\ell+1}= (\sqrt{-1})^{-k-1},\\
&& b_j=1,
j=1,\ldots,2\ell+1,\quad X=-2\sqrt{-1}\pi k;
\end{eqnarray*}
$\Phi\circ \tilde{\delta}: S^1 \to \ymS{2\ell+1}{0}_k$ is a loop
defined by
$$
a_j=b_j=1, j=1, \ldots,2\ell,$$
$$ a_{2\ell+1}=(\sqrt{-1})^{-k-1}
e^{\sqrt{-1}\theta},\quad b_{2\ell+1}=1, \quad X=-2\sqrt{-1}\pi k;
$$
and $\Phi\circ \tilde{\gamma}: S^1 \to \ymS{2\ell+1}{0}_k$ is a loop
defined by
$$
a_j=b_j=1, j=1,\ldots,2\ell,
$$
$$ a_{2\ell+1}=(\sqrt{-1})^{-k-1}, \quad
b_{2\ell+1}=e^{\sqrt{-1}\theta},\quad X=-2\sqrt{-1}\pi k.
$$
Thus the $(4\ell+2)$ loops
\begin{equation}\label{eqn:oneKloop}
\Phi\circ \tilde{\alpha}_i,\ \Phi\circ \tilde{\beta}_i,\ \Phi\circ
\tilde{\alpha}'_i,\ \Phi\circ \tilde{\beta}'_i,\
i=1,\ldots,\ell,\quad \Phi \circ \tilde{\delta},\  \Phi\circ
\tilde{\gamma}
\end{equation}
generate the fundamental group of $\ymS{2\ell+1}{0}_k\cong
U(1)^{4\ell+2}$.

To show that $V_{1,k} \to \ymS{2\ell+1}{0}_k$ is orientable, it
suffices to show that its restriction to each of the $(4\ell+2)$
loops in \eqref{eqn:oneKloop} is orientable, or equivalently:
\begin{pro}\label{thm:oneK}
Let $\ell\geq 0$. The restriction of $W_{1,k}\to \ZymS{2}_k$ to each
of the following $(4\ell+2)$ loops is orientable:
$$
\tilde{\alpha}_i,\ \tilde{\beta}_i,\ \tilde{\alpha}'_i,\
\tilde{\beta}'_i,\ i=1,\ldots,\ell,\quad \tilde{\delta},\
\tilde{\gamma}.
$$
\end{pro}
\begin{proof}
See Section \ref{sec:Kproof}.
\end{proof}

\subsubsection{$n>1$, $i=1$} \label{sec:nRPloop}
Let $\ell\geq 1$. Let $\alpha_i,\ \beta_i,\ \alpha'_i,\ \beta'_i$
be defined in Section \ref{sec:RPloop}, so that the holonomies along
$\alpha_i,\ \beta_i,\ \alpha'_i,\ \beta'_i$ are $a_i,\ b_i,\ a'_i,\
b'_i\in U(n)$, respectively.

Let $a_\theta =
\diag(e^{\sqrt{-1}\theta},\underbrace{1,\ldots,1}_{n-1})\in U(n)$.
By Goto's commutator theorem, the map $G^2 \to G$ defined by
$(a,b)\mapsto [a,b]$ is surjective if $G$ is semisimple (cf. \cite[Theorem 9.2]{HM}). So there
exist $a,b\in SU(n)$ such that
$$
[a,b]=e^{2\pi\sqrt{-1}k/n}I_n \in SU(n).
$$
Let $\tilde{\alpha}_i :S^1\to \ZymU{1}_{\kn}$ be the loop defined by
\begin{eqnarray*}
&&a_i=a_\theta,\quad b_i=I_n,\quad a_i'= a,\quad b_i'=b,\\
&& a_j=b_j=a_j'=b_j'=I_n \textup{ for }j\neq i, \quad
c=c'=e^{\frac{\pi\sqrt{-1}k}{2n}} I_n, \quad
X=-2\sqrt{-1}\pi\frac{k}{n}I_n.
\end{eqnarray*}
We define $\tilde{\beta}_i$, $\tilde{\alpha}'_i$, $\tilde{\beta}'_i$
similarly.  From \eqref{eqn:PhiRP} it is
clear that the fundamental group of $\ymS{2\ell}{0}_k \cong
U(1)^{4\ell}$ is generated by the following $4\ell$ loops:
$$
\Psi\circ \Phi \circ \tilde{\alpha}_i,\
\Psi\circ\Phi\circ\tilde{\beta}_i,\ \Psi\circ \Phi\circ
\tilde{\alpha}'_i,\ \Psi\circ \Phi \circ \tilde{\beta}'_i,\quad
i=1,\ldots,\ell.
$$
By Proposition \ref{thm:Xdet}, the fundamental group of
$\ymU{2\ell}{0}_\kn$ is generated by the following $4\ell$ loops:
\begin{equation} \label{eqn:nRPloop}
\Phi \circ \tilde{\alpha}_i,\ \Phi\circ\tilde{\beta}_i,\ \Phi\circ
\tilde{\alpha}'_i,\ \Phi \circ \tilde{\beta}'_i,\quad
i=1,\ldots,\ell.
\end{equation}
To show that $V_{n,k} \to \ymU{2\ell}{0}_\kn$ is orientable, it
suffices to show that its restriction to each of the $4\ell$ loops
in \eqref{eqn:nRPloop} is orientable, or equivalently:
\begin{pro}\label{thm:nRP}
Let $\ell\geq 1$ and let $n>1$. The restriction of the vector bundle $W_{n,k}\to
\ZymU{1}_\kn$ to each of the following $4\ell$ loops is orientable:
$$
\tilde{\alpha}_i,\ \tilde{\beta}_i,\ \tilde{\alpha}'_i,\
\tilde{\beta}'_i,\quad i=1,\ldots,\ell.
$$
\end{pro}
\begin{proof}
See Section \ref{sec:RPproof}.
\end{proof}

\subsubsection{$n>1$, $i=2$}
Let $\ell\geq 1$. Let $\alpha_i,\ \beta_i,\ \alpha'_i,\ \beta'_i,\
\delta,\ \delta',\ \gamma,\ \gamma'$ be defined as in Section
\ref{sec:Kloop}, so that the holonomies along $\alpha_i,\ \beta_i,\
\alpha'_i,\ \beta'_i,\ \delta,\ \delta',\ \gamma,\ \gamma'$ are
$a_i,\ b_i,\ a'_i,\ b'_i,\ d,\ d',\ c,\ c' \in U(n)$, respectively.
Define $a,b,a_\theta$ as in Section \ref{sec:nRPloop}, so that
$[a,b]=e^{2\sqrt{-1}\pi k/n}I_n$.

Let $\tilde{\alpha}_i: S^1\to \ZymU{2}_\kn$ be the loop defined by
\begin{eqnarray*}
&&
a_i=a_\theta,\quad b_i = I_n, \quad  a_i'=a,\quad b_i'=b,\\
&& a_j=b_j=a_j'=b_j'=I_n \textup{ for }j\neq i, \quad c=c'=I_n,\\
&& d=\sqrt{-1}e^{\frac{\pi\sqrt{-1}k}{2n}} I_n,
d'=-\sqrt{-1}e^{\frac{\pi\sqrt{-1}k}{2n}} I_n, \quad
X=-2\sqrt{-1}\pi\frac{k}{n}I_n.
\end{eqnarray*}
We $\tilde{\beta}i$, $\tilde{\alpha}'_i$, $\tilde{\beta}'_i$
similarly.

Let $\tilde{\delta}:S^1\to \ZymU{2}_\kn$ be the loop defined by
\begin{eqnarray*}
&& a_1=b_1=I_n,\quad a_1'=a,\quad b_1'=b,\quad d= \sqrt{-1}
e^{\frac{\pi\sqrt{-1}k}{2n} }\bar{a}_\theta,
d'=-\sqrt{-1}e^{\frac{\pi\sqrt{-1}k}{2n}} a_\theta,\\
&&  c=c'=I_n,\quad  a_j=b_j=a_j'=b_j'=I_n,\ j=2,\ldots,\ell,\quad
X=-2\sqrt{-1}\pi\frac{k}{n}I_n.
\end{eqnarray*}

Let $\tilde{\gamma}: S^1\to \ZymU{2}_\kn$ be the loop defined by
\begin{eqnarray*}
&& a_1=b_1=I_n,\quad a_1'=a, \quad b_1'=b,\quad c=a_\theta,\quad c'=I_n,\\
&& d=\sqrt{-1}e^{\frac{\pi\sqrt{-1}k}{2n} }I_n,
\quad d'=-\sqrt{-1} e^{\frac{\pi\sqrt{-1}k}{2n}}I_n,\\
&& a_j=b_j=a_j'=b_j'=I_n, j=2,\ldots,\ell,\quad
X=-2\sqrt{-1}\pi\frac{k}{n}I_n.
\end{eqnarray*}

From \eqref{eqn:PhiK}, it is clear that the following $4\ell+2$ loops generate the fundamental group
of $\ymS{2\ell+1}{0}_k \cong U(1)^{4\ell+2}$:
$$
\Psi\circ\Phi\circ \tilde{\alpha}_i,\ \Psi\circ\Phi\circ
\tilde{\beta}_i,\ \Psi\circ\Phi\circ \tilde{\alpha}'_i,\
\Psi\circ\Phi\circ \tilde{\beta}'_i,\ i=1,\ldots,\ell,\quad
\Psi\circ \Phi \circ \tilde{\delta},\  \Psi\circ\Phi\circ
\tilde{\gamma}.
$$

By Proposition \ref{thm:Xdet}, the fundamental group of
$\ymU{2\ell+1}{0}_\kn$ is generated by the following $4\ell+2$
loops:
\begin{equation}\label{eqn:nKloop}
\Phi\circ \tilde{\alpha}_i,\ \Phi\circ \tilde{\beta}_i,\ \Phi\circ
\tilde{\alpha}'_i,\ \Phi\circ \tilde{\beta}'_i,\
i=1,\ldots,\ell,\quad \Phi \circ \tilde{\delta},\  \Phi\circ
\tilde{\gamma}.
\end{equation}

To show that $V_{n,k}\to \ymU{2\ell+1}{0}_\kn$ is orientable, it
suffices to show that its restriction to each of the $(4\ell+2)$
loops in \eqref{eqn:nKloop} is orientable, or equivalently:
\begin{pro}\label{thm:nK}
Let $\ell\geq 1$ and let $n>1$. The restriction of the vector bundle
$W_{n,k}\to \ZymU{2}_\kn$ to each of the following $(4\ell+2)$ loops is
orientable:
$$
\tilde{\alpha}_i,\ \tilde{\beta}_i,\ \tilde{\alpha}'_i,\
\tilde{\beta}'_i,\ i=1,\ldots,\ell,\quad \tilde{\delta},\
\tilde{\gamma}.
$$
\end{pro}
\begin{proof}
See Section \ref{sec:Kproof}.
\end{proof}

\section{Orientability along Loops} \label{sec:alongloops}

Our approach is similar to that in the proof of \cite[Proposition
21.3]{FOOO}.  Let $\gamma : S^1 \to \ZymU{i}_{\kn}$ be any of the loops in Proposition
\ref{thm:oneRP}---\ref{thm:nK}. We need to show that $E:= \gamma^*W_{n,k}$
is an orientable real vector bundle over $S^1$.  The fiber of $E$ over $\theta\in S^1$
is given by $E_\theta = H^1(\tSi,\cM_\theta)^{\tau_\theta}$, where
$\cM_\theta$ is an $S^1$-family of holomorphic vector bundles
over $\tSi$. Our strategy is to use degeneration and
normalization of the Riemann surface $\tSi$ to
show that $E\cong E_\bC \oplus  E_\bR$, where
$E_\bC$ is a complex vector bundle over $S^1$, and $E_\bR$ is
a (possibly zero, possibly virtual) real vector bundle
which we can describe very explicitly. The explicit description of
$E_\bR$ allows us to compute $w_1(E_\bR)=0$.

\subsection{Degeneration of the Riemann surface}\label{sec:degeneration}
We degenerate the smooth Riemann surface
$\tilde{\Si}=\Si_0^{2\ell+i-1}$ to a nodal Riemann surface $C$ with
three irreducible components $C_+$, $C_0$ and $C_-$, where
$C_\pm\cong \Si^\ell_0$, $C_0\cong \Si^0_i$, and $C_0$ intersects
$C_\pm$ at a node $p_\pm$. More precisely, consider a family of
Riemann surfaces $\tilde{\Si}_t$, where $t\in I= [0,1]$, such that
\begin{enumerate}
\item[(i)] $\tSi_1 = \tSi$, $\tSi_0=C$.
\item[(ii)] $\tSi_t$ is smooth for $0< t\leq 1$.
\item[(iii)] There is a family of anti-holomorphic involutions
$\tau_t: \tSi_t \to \tSi_t$, such that
$$
\tau_0(C_\pm)= C_\mp,\quad \tau_0(p_\pm)= p_\mp,\quad
\tau_0(C_0)=C_0,\quad \tau_1 =\tau.
$$
\end{enumerate}
The $i=1$ and $i=2$ cases are shown in Figure 1 and Figure 2,
respectively. In Figure 1 and 2, $\tau_t(q_\pm)=q_\mp$,
$\tau_t(r_\pm)=r_\mp$, $\tau_t(\ep_\pm)=\ep_\mp$. Notice that our
loops start from $p_{\pm}$, so the loop $\alpha_1$ (resp.
$\alpha_1'$) contains the path from $p_+$ to $q_+$ (resp. from $p_-$
to $q_-$) and its inverse; the loop $\beta_1$ (resp. $\beta_1'$)
contains the path from $p_+$ to $r_+$ (resp. from $p_-$ to $r_-$)
and its inverse. In the degeneration $t\to 0$, the loop $\ep_\pm$
shrinks to the point $p_\pm$, respectively.

\begin{figure}[h]
\begin{center}
\psfrag{p+}{\tiny $p_+$} \psfrag{p-}{\tiny $p_-$} \psfrag{q+}{\tiny
$q_+$} \psfrag{q-}{\tiny $q_-$} \psfrag{r+}{\tiny $r_+$}
\psfrag{r-}{\tiny $r_-$} \psfrag{a}{\tiny $\beta_1$}
\psfrag{b}{\tiny $\alpha_1$} \psfrag{ap}{\tiny $\beta_1'$}
\psfrag{bp}{\tiny $\alpha_1'$} \psfrag{r}{\tiny $\gamma$}
\psfrag{rp}{\tiny $\gamma'$} \psfrag{S0}{\small $\tSi_0$}
\psfrag{St}{\small $\tSi_t$} \psfrag{t>0}{\tiny $(t>0)$}
\psfrag{e+}{\tiny $\ep_+$} \psfrag{e-}{\tiny $\ep_-$}
\psfrag{C0}{\small $C_0$} \psfrag{C+}{\small $C_+$}
\psfrag{C-}{\small $C_-$}
\includegraphics[scale=0.8]{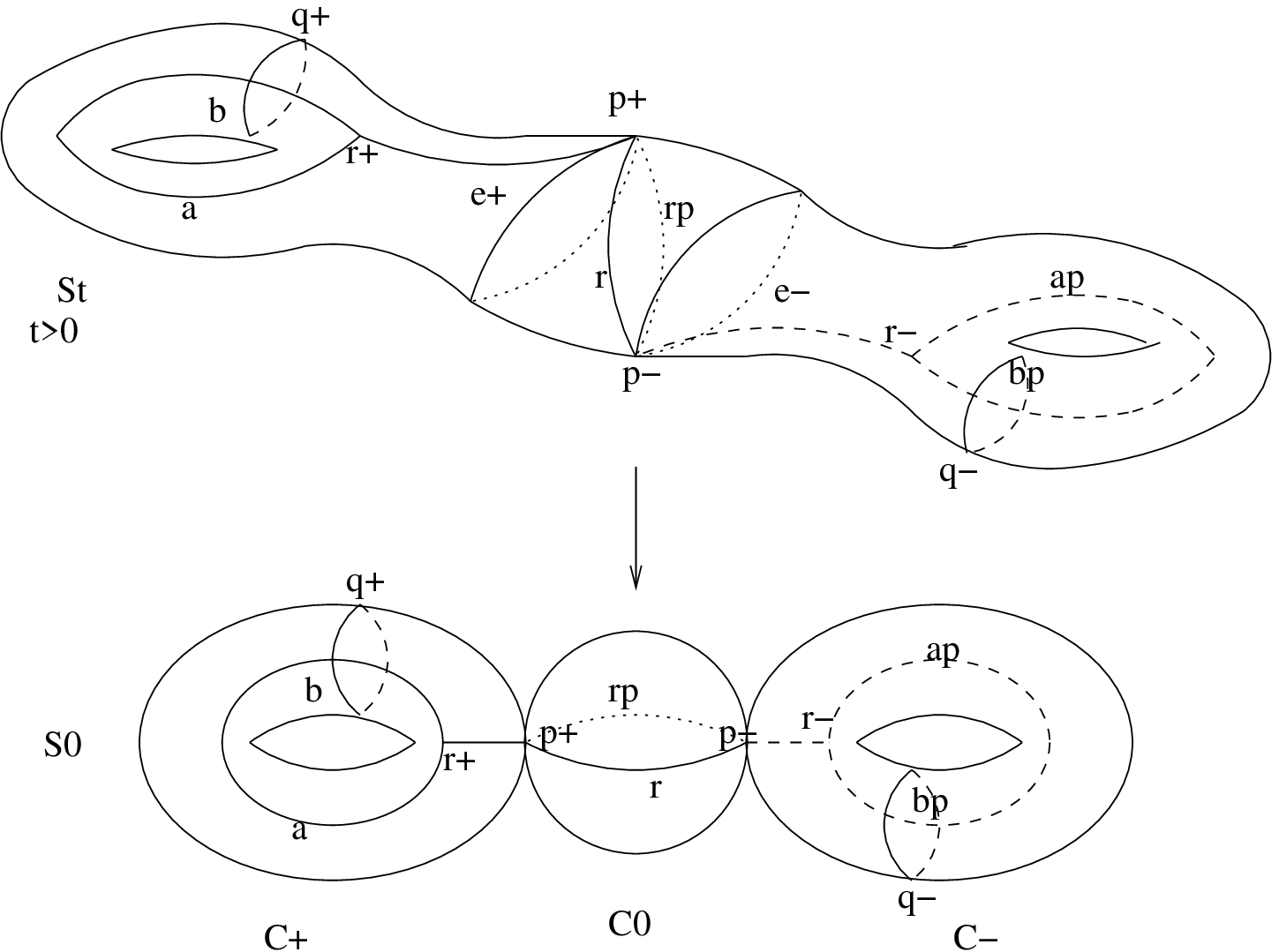}
\end{center}
\caption{Degeneration of $\widetilde{\Si^\ell_1}= \Sigma^{2\ell}_0$.
}
\end{figure}

\begin{figure}[h]
\begin{center}
\psfrag{p+}{\tiny $p_+$} \psfrag{p-}{\tiny $p_-$} \psfrag{q+}{\tiny
$q_+$} \psfrag{q-}{\tiny $q_-$} \psfrag{r+}{\tiny $r_+$}
\psfrag{r-}{\tiny $r_-$} \psfrag{a}{\tiny $\beta_1$}
\psfrag{b}{\tiny $\alpha_1$} \psfrag{ap}{\tiny $\beta_1'$}
\psfrag{bp}{\tiny $\alpha_1'$} \psfrag{r}{\tiny $\gamma$}
\psfrag{rp}{\tiny $\gamma'$} \psfrag{d}{\tiny $\delta$}
\psfrag{dp}{\tiny $\delta'$} \psfrag{S0}{\small $\tSi_0$}
\psfrag{St}{\small $\tSi_t$} \psfrag{t>0}{\tiny $(t>0)$}
\psfrag{e+}{\tiny $\ep_+$} \psfrag{e-}{\tiny $\ep_-$}
\psfrag{C0}{\small $C_0$} \psfrag{C+}{\small $C_+$}
\psfrag{C-}{\small $C_-$}
\includegraphics[scale=0.8]{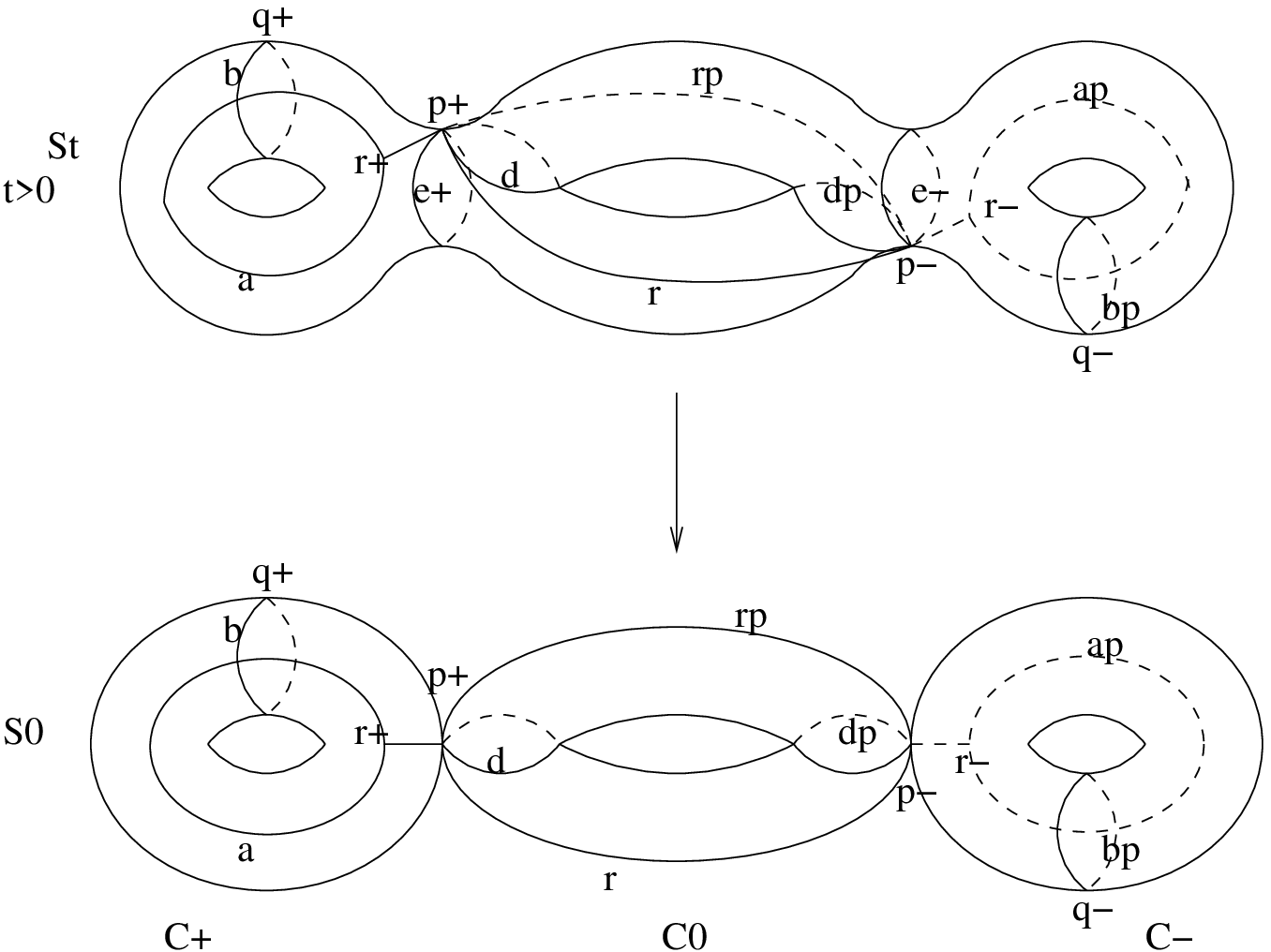}
\end{center}
\caption{Degeneration of $\widetilde{\Si^\ell_2}=
\Sigma^{2\ell+1}_0$. }
\end{figure}

In Section \ref{sec:RPproof} and Section \ref{sec:Kproof}
we will degenerate the family $\cM_\theta$ together
with the base $\tSi$ to obtain an $S^1$-family of vector bundles
$\cM_{\theta,0}\to \tSi_0=C$. Then we will reduce the orientability
of $E\to S^1$ to the orientability of
$E_0\to S^1$ whose fiber over $\theta\in S^1$
is $H^1(C,\cM_{\theta,0})^{\tau_{\theta,0}}$.

\subsection{Normalization}

Let $C_+$, $C_0$, $C_-$, and $C$ be defined as in Section \ref{sec:degeneration}.
The normalization $\tilde{C}$ of $C$ is a disconnected smooth Riemann surface
which can be identified with the disjoint union of $C_+$, $C_0$, and $C_-$.
There is a normalization map $\tilde{\nu}:\tilde{C}\to C$, identifying a point on $C_\pm$
to a point in $C_0$ (which becomes the node $p_\pm$).
We will use the following lemma to study the orientability of $E_0\to S^1$.

\begin{lm}\label{lm:LES} Let  $\cM \to C$ be a holomorphic vector bundle of rank $r$.
Let $\cM_+ \to C_+$, $\cM_0\to C_0$ and $\cM_- \to C_-$ be the restrictions
of $\cM$ to $C_+$, $C_0$, and $C_-$, respectively. Then we have
a long exact sequence of complex vector spaces:
\begin{equation}\label{eqn:LES}
\begin{aligned}
0 &\to H^0(C,\cM)\to H^0(C_+,\cM_+)\oplus H^0(C_0,\cM_0)\oplus H^0(C_-,\cM_-)\\
& \to \bC^r_{p_+}\oplus \bC^r_{p_-}\to H^1(C,\cM)\\
& \to H^1(C_+,\cM_+)\oplus H^1(C_0,\cM_0)\oplus H^1(C_-,\cM_-) \to 0.
\end{aligned}
\end{equation}
\end{lm}
\begin{proof}  Let
$\iota_+:C_+\injects C$, $\iota_0:C_0\injects C$, and
$\iota_-:C_-\injects C$ be inclusion maps, and let $\nu:\tilde{C}\to C$
be the normalization map. Then
$$
\nu_* \cO_{\tilde{C}}= \iota_{+*} \cO_{C_+}\oplus \iota_{0*}\cO_{C_0} \oplus \iota_{-*}\cO_{C_-}.
$$
We have a short exact sequence of sheaves on $C$ (known as the
normalization sequence):
\begin{equation}\label{eqn:normalization}
0\to \cO_C\to \iota_{+*}\cO_{C_+} \oplus \iota_{0*}\cO_{C_0}\oplus \iota_{-*}\cO_{C-}
\to \cO_{p_+}\oplus \cO_{p_-}\to 0.
\end{equation}
(See e.g. page 81 of \cite{HaM}.) Twisting the normalization sequence
\eqref{eqn:normalization}  by $\cM$, we obtain a short exact sequence of
sheaves on $C$:
\begin{equation}\label{eqn:SES}
0\to \cM\to \iota_{+*}\cM_+ \oplus \iota_{0*}\cM_0\oplus \iota_{-*}\cM_-
\to \cO^{\oplus r}_{p_+}\oplus \cO^{\oplus r}_{p_-}\to 0.
\end{equation}
The long exact sequence of cohomology groups associated to \eqref{eqn:SES} is
\begin{equation}\label{eqn:LES-I}
\begin{aligned}
0 &\to H^0(C,\cM)\to H^0(C,\iota_{+*}\cM_+)\oplus H^0(C,\iota_{0*}\cM_0)\oplus H^0(C,\iota_{-*}\cM_-)\\
& \to \bC^r_{p_+}\oplus \bC^r_{p_-}\to H^1(C,\cM)\\
& \to H^1(C,\iota_{+*}\cM_+)\oplus H^1(C,\iota_{0*}\cM_0)\oplus H^1(C,\iota_{-*}\cM_-) \to 0.
\end{aligned}
\end{equation}
For $k=0,1$, we have
$$
H^k(C, \iota_{\pm *}\cM_{\pm}) = H^k(C_\pm, \cM_\pm),\quad
H^k(C, \iota_{0 *}\cM_0)=H^k(C_0,\cM_0).
$$
So \eqref{eqn:LES-I} is equivalent to \eqref{eqn:LES}.
\end{proof}

\subsection{The $i=1$ case}\label{sec:RPproof}
\begin{proof}[Proof of Proposition \ref{thm:oneRP} and Proposition \ref{thm:nRP}]
We will show that $\tilde{\alpha}_j^* W_{n,k} \to S^1$ is
orientable, $j=1,\ldots,\ell$. The other loops $\tilde{\beta}_j,
\tilde{\alpha}'_j,\tilde{\beta}'_j$ are similar.

Note that $\tilde{\alpha}_j^* W_{n,k} = (\phi_T\circ
\tilde{\alpha}_j)^* U_{n,k}^{\hat{\tau}}$. We have
$$
\tilde{\alpha}_j(\theta) = (V_\theta, e^{\sqrt{-1}\pi k/2n} I_{n},
V', e^{\sqrt{-1}\pi k/2n} I_n, -2\pi\sqrt{-1}\frac{k}{n} I_n),
$$
where $V_\theta, V' \in U(n)^{2\ell}$. Note that this is also true
for the $n=1$ case.  The loop $\phi_T\circ \tilde{\alpha}_j: S^1 \to
Z_{\mathrm{YM}}^{\ell,1}(U(n^2))_{\mtwokn}$ is given by (see
\eqref{eqn:RPphiT} for the definition of $\phi_T$):
\begin{equation}\label{eqn:MholRP}
\phi_T\circ \tilde{\alpha}_j(\theta) =(\bV_\theta\otimes V', T,
V_\theta \otimes \bV', T, 4\pi\sqrt{-1}\frac{k}{n} I_{n^2}).
\end{equation}

The loop $\phi_T\circ \tilde{\alpha}_j(\theta)$ can be viewed as an
$S^1$-family of polystable holomorphic vector bundles $\cM_\theta$
of rank $n^2$, degree $-2kn$ over the Riemann surface $\tSi\cong
\Si^{2\ell}_0$. We now consider holomorphic vector bundles
$\cM_{\theta,t}$ of rank $n^2$, degree $-2kn$ over $\tSi_t$ with the
following properties:
\begin{enumerate}
\item[(i)] For $0< t\leq 1$, the holonomies of $\cM_{\theta, t}$
over $\tSi_t$ are given by \eqref{eqn:MholRP}.
\item[(ii)] When $t=0$, we have
$$
\cM_{\theta,0}\Bigr|_{C_\pm} = \cM_\pm(\theta),\quad
\cM_{\theta,0}\Bigr|_{C_0}= \cM_0(\theta),
$$
where $\cM_\pm(\theta)$ is a rank $n^2$, degree $-kn$ polystable
holomorphic vector bundle over $C_\pm\cong \Si^\ell_0$ and
$\cM_0(\theta)$ is a rank $n^2$, degree $0$ polystable holomorphic
vector bundle over $C_0\cong\bP^1$.

\item[(iii)] The holonomies of $\cM_+(\theta)$ starting from $p_+$
 along $(\alpha_1,\beta_1,...,\alpha_\ell,\beta_\ell)$ are given by
$$
(\bV_\theta \otimes V', 2\pi\sqrt{-1}\frac{k}{n}I_{n^2}) \in
X_{\mathrm{YM}}^{\ell,0}(U(n^2))_{\mkn},
$$
and the holonomies of $\cM_-(\theta)$ starting from $p_-$ along
$(\beta'_\ell,\alpha'_\ell,..,\beta'_1,\alpha'_1)$ is given by $$
(\fr(V_\theta\otimes \bV'),2\pi\sqrt{-1}\frac{k}{n}I_{n^2}) \in
X_{\mathrm{YM}}^{\ell,0}(U(n^2))_{\mkn}.$$

\item[(iv)] By (ii),  $\cM_0(\theta) \cong \cO_{\bP^1}^{\oplus n^2}$, since
the trivial bundle is the only degree zero polystable bundle on
$\bP^1$. Thus, $\cM_0(\theta)$ is independent of $\theta$ and will
be denoted by $\cM_0$. The holonomy of $\cM_0$ along the equator
$\gamma\gamma'$ is $TT=I_{n^2}$ as expected.

\item[(v)] For all $(\theta,t)\in S^1\times I$,
we have $\overline{\tau_t^*\cM_{\theta,t}} = \cM_{\theta, t}$, so
there is a conjugate linear involution $\hat{\tau}_{\theta,t}$ on
$H^*(\tSi_t, \cM_{\theta,t})$.
\end{enumerate}

As $(\theta,t)$ varies, the real vector spaces
$$
\{ H^1(\tSi_t,\cM_{\theta,t})^{\hat{\tau}_{\theta,t}} \mid
(\theta,t)\in S^1 \times I\}
$$
form a real vector bundle $E$ over the cylinder $S^1 \times I$. Let
$i_t: S^1 \to S^1\times I$ be the embedding $\theta\mapsto
(\theta,t)$. Then
$$
i_1^* E = (\phi_T\circ\tilde{\alpha}_j)^* U_{n,k}^{\hat{\tau}} =
\tilde{\alpha}_j^* W_{n,k}.
$$
The maps $i_0$ and $i_1$ homotopic, so $i_0^*E$ and $i_1^*E$ are
isomorphic real vector bundles over $S^1$.  Thus
$\tilde{\alpha}_j^*W_{n,k}$ is orientable if and only if $i_0^* E$
is an orientable real vector bundle over $S^1$.

By Lemma \ref{lm:LES}, we have the following long exact sequence:
\begin{equation}\label{eqn:RP-long}
\begin{aligned}
0 \to &H^0(C,\cM_{\theta,0})\to H^0(C_+,\cM_+(\theta))\oplus
H^0(C_0,\cM_0)\oplus H^0(C_-,\cM_-(\theta))\\
&\to \bC^{n^2}_{p_+}\oplus \bC^{n^2}_{p_-}
\to H^1(C,\cM_{\theta,0})\\
& \to H^1(C_+,\cM_+(\theta))\oplus
H^1(C_0,\cM_0)\oplus H^1(C_-,\cM_-(\theta))\to 0
\end{aligned}
\end{equation}
where $p_\pm$ is the node at which $C_0$ and $C_\pm$ intersect. We
have
$$
\deg\cM_{\theta,0}=-2nk <0,\quad \deg\cM_\pm(\theta)=-nk<0,
$$
so
$$
H^0(C,\cM_{\theta,0})=H^0(C_+,\cM_+(\theta))=H^0(C_-,\cM_-(\theta))=0.
$$

By (iv), $\cM_0\cong \cO^{\oplus n^2}_{\bP^1}$, so $H^0(C_0,\cM_0)\cong
\bC^{\oplus n^2}$ and $H^1(C_0,\cM_0)=0$. The map $H^0(C_0,\cM_0)\to
\bC^{n^2}_{p_\pm}$ is the evaluation map $s\mapsto s(p_\pm)$.
Therefore \eqref{eqn:RP-long} is reduced to
\begin{eqnarray}\label{eqn:RP-short}
0 \to \bC^{n^2} \stackrel{j}{\to} \bC_{p_+}^{\oplus n^2} \oplus
\bC_{p_-}^{\oplus n^2}\to H^1(C,\cM_{\theta,0})\to\notag\\
H^1(C_+,\cM_+(\theta))\oplus H^1(C_-,\cM_-(\theta))\to 0
\end{eqnarray}
where $j(v)=(v,v)$. The involution $\hat{\tau}_{\theta,0}$ acts on
the exact sequence \eqref{eqn:RP-short} in the following way :
$$
H^1(C,\cM_{\theta,0})
\stackrel{\hat{\tau}_{\theta,0}}{\longleftrightarrow}
H^1(C,\cM_{\theta,0}),\quad
H^1(C_+,\cM_+(\theta))\stackrel{\hat{\tau}_{\theta,0}}{\longleftrightarrow}
H^1(C_-,\cM_-(\theta)).
$$

The involution on $\bC^{\oplus n^2}$ and  $\bC^{\oplus
n^2}_{p_+}\oplus \bC^{\oplus n^2}_{p_-}$ is independent of $\theta$:
\begin{eqnarray*}
&& \hat{\tau}_{\theta,0}:\bC^{n^2}\to \bC^{n^2}, \quad v\mapsto \bar{v}\\
&& \hat{\tau}_{\theta,0}:\bC_{p_+}^{n^2}\oplus \bC^{n^2}_{p_-}\to
\bC_{p_+}^{n^2}\oplus\bC_{p_-}^{n^2}, \quad (v,w)\mapsto
(\bar{w},\bar{v}).
\end{eqnarray*}
Thus, we have
$$
0\to \bR^{n^2}\to H^1(C,\cM_{\theta,0})^{\hat{\tau}_\theta,0} \to
H^1(C_+,\cM_+(\theta))\to 0.
$$

We conclude that $i_0^*E \cong E_\bC \oplus E_\bR$, where $E_\bC\to
S^1$ is a complex vector bundle whose fiber at $\theta\in S^1$ is
$H^1(C_+,\cM_+(\theta))$, and $E_\bR\to S^1$ is a trivial real
vector bundle of rank $n^2$. Therefore $i_0^*E$ is orientable.
\end{proof}

\subsection{The $i=2$ case}\label{sec:Kproof}
\begin{proof}[Proof of Proposition \ref{thm:oneK} and Proposition \ref{thm:nK}]
We will study the orientability of the real vector bundle
$$
\lambda^*W_{n,k} = (\phi_T\circ\lambda)^* U_{n,k}^{\hat{\tau}}
$$
over $S^1$, where $\lambda: S^1 \to \ZymU{2}_{\kn}$ is one of the
$4\ell +2 $ loops. Note that $\lambda$ is of the form
$$
\lambda(\theta) = (V_\theta, \sqrt{-1}e^{\frac{\pi\sqrt{-1}k}{2n} }
g_\theta ,c_\theta, V', -\sqrt{-1}e^{\frac{\pi\sqrt{-1}k}{2n}}
\bar{g}_\theta, I_n, -2\pi\sqrt{-1}\frac{k}{n} I_n),
$$
where $V_\theta, V'_\theta \in U(n)^{2\ell}$ and $g_\theta$,
$c_\theta$ are diagonal matrices in $U(n)$. Note that this is also
true when $n=1$.

The loop $\phi_T\circ \lambda: S^1 \to
Z_{\mathrm{YM}}^{\ell,1}(U(n^2))_{\mtwokn}$ is given by (see
\eqref{eqn:KphiT} for the definition of $\phi_T$):
\begin{equation}\label{eqn:MholK}
\begin{aligned}
\phi_T\circ \lambda(\theta) =& (\bV_\theta\otimes V',
-\bar{g}_\theta \otimes \bar{g}_\theta,
\left(\bc_\theta\otimes I_n \right)T, \\
& V_\theta\otimes \bV', - g_\theta\otimes g_\theta,
\left(c_\theta\otimes I_n \right) T, 4\pi\sqrt{-1}\frac{k}{n}
I_{n^2}).
\end{aligned}
\end{equation}

The loop $\phi_T\circ \lambda(\theta)$ can be viewed as an
$S^1$-family of polystable holomorphic vector bundles $\cM_\theta$
of rank $n^2$, degree $-2kn$ over the Riemann surface $\tSi\cong
\Si^{2\ell+1}_0$. We now consider holomorphic vector bundles
$\cM_{\theta,t}$ of rank $n^2$, degree $-2kn$ over $\tSi_t$ with the
following properties:
\begin{enumerate}
\item[(i)] For $0< t\leq 1$, the holonomies of $\cM_{\theta, t}$
are given by \eqref{eqn:MholK}.
\item[(ii)] When $t=0$, we have
$$
\cM_{\theta,0}\Bigr|_{C_\pm} = \cM_\pm(\theta),\quad
\cM_{\theta,0}\Bigr|_{C_0}= \cM_0(\theta),
$$
where $\cM_\pm(\theta)$ is a rank $n^2$, degree $-kn$ polystable
holomorphic vector bundle over $C_\pm\cong \Si^\ell_0$ and
$\cM_0(\theta)$ is a rank $n^2$, degree $0$ polystable holomorphic
vector bundle over $C_0\cong S^1\times S^1$.

\item[(iii)]
The holonomies of $\cM_+(\theta)$ starting from $p_+$
 along $(\alpha_1,\beta_1,...,\alpha_\ell,\beta_\ell)$ are given by
$$
(\bV_\theta \otimes V', 2\pi\sqrt{-1}\frac{k}{n}I_{n^2}) \in
X_{\mathrm{YM}}^{\ell,0}(U(n^2))_{\mkn},
$$
and the holonomies of $\cM_-(\theta)$  along
$(\beta'_\ell,\alpha'_\ell,..,\beta'_1,\alpha'_1)$ (starting from
$p_-$) are given by $$
(\fr(V_\theta\otimes \bV'),2\pi\sqrt{-1}\frac{k}{n}I_{n^2}) \in
X_{\mathrm{YM}}^{\ell,0}(U(n^2))_{\mkn}.$$

\item[(iv)] The holonomies of $\cM_0(\theta)$ along
$\delta$, $\gamma$, $\delta'$, $\gamma'$ are given by
\begin{equation}\label{eqn:lambda-zero}
\lambda_0(\theta)= (-\bar{g}_\theta\otimes \bar{g}_\theta,
(\bc_\theta\otimes I_n) T, -g_\theta\otimes g_\theta,
(c_\theta\otimes I_n)T ) \in Z_{\mathrm{flat}}^{0,2}(U(n^2))^{\hat{\tau}}.
\end{equation}
Therefore, the holomomies of $\cM_0(\theta)$ along
the loops $\delta^{-1}$, $\gamma\gamma'$ are given by
$$
\Phi^{0,2}\circ \lambda_0(\theta)= (-g_\theta\otimes g_\theta, \bc_\theta\otimes c_\theta)
\in X_{\mathrm{flat}}^{1,0}(U(n^2)).
$$

\item[(v)] For all $(\theta,t)\in S^1\times I$,
we have $\overline{\tau_t^*\cM_{\theta,t}} = \cM_{\theta, t}$, so
there is a conjugate linear involution $\hat{\tau}_{\theta,t}$ on
$H^*(\tSi_t, \cM_{\theta,t})$.
\end{enumerate}

As $(\theta,t)$ varies, the real vector spaces
$$
\{ H^1(\tSi_t,\cM_{\theta,t})^{\hat{\tau}_{\theta,t}} \mid
(\theta,t)\in S^1 \times I\}
$$
form a real vector bundle $E$ over the cylinder $S^1 \times I$. Let
$i_t: S^1 \to S^1\times I$ be the embedding $\theta\mapsto
(\theta,t)$. Then
$$
i_1^* E = (\phi_T\circ\lambda)^* U_{n,k}^{\hat{\tau}} = \lambda^*
W_{n,k}.
$$
The maps $i_0$ and $i_1$ are homotopic, so $i_0^*E$ and $i_1^*E$ are
isomorphic real vector bundles over $S^1$.  Thus $\lambda^*W_{n,k}$
is orientable if and only if $i_0^* E$ is an orientable vector
bundle over $S^1$.

By Lemma \ref{lm:LES}, we have the following long exact sequence:
\begin{equation}\label{eqn:K-long}
\begin{aligned}
0 \to &H^0(C,\cM_{\theta,0})\to H^0(C_+,\cM_+(\theta))\oplus
H^0(C_0,\cM_0(\theta))\oplus H^0(C_-,\cM_-(\theta))\\
&\to \bC_{p_+}^{n^2}\oplus \bC_{p_-}^{n^2}
\to H^1(C,\cM_{\theta,0})\\
&\to H^1(C_+,\cM_+(\theta))\oplus
H^1(C_0,\cM_0(\theta))\oplus H^1(C_-,\cM_-(\theta))\to 0
\end{aligned}
\end{equation}
where $p_\pm$ is the node at which $C_0$ and $C_\pm$ intersect. We
have
$$
\deg\cM_{\theta,0}=-2nk <0,\quad \deg\cM_\pm(\theta)=-nk<0,
$$
so
$$
H^0(C,\cM_{\theta,0})=H^0(C_+,\cM_+(\theta))=H^0(C_-,\cM_-(\theta))=0.
$$

Note that the holonomies of $\cM_0(\theta)$ are diagonal, so it is
the direct sum of $n^2$ holomorphic line bundles of degree $0$. Let
$\LL_{a,b}$ denote the degree $0$ holomorphic line bundle whose
holonomies along the loops $\delta^{-1}$, $\gamma\gamma'$ are given by
$$
(a,b)\in \flS{1}{0}=U(1)^2.
$$
Then
$$
H^0(C_0,\LL_{a,b})=H^1(C_0,\LL_{a,b})
=\begin{cases}
\bC,\quad (a,b)=(1,1)\\
0,\quad (a,b)\neq (1,1).
\end{cases}
$$

\begin{enumerate}
\item[Case 1.]
$\lambda=\tilde{\alpha}_i, \tilde{\beta}_i, \tilde{\alpha}'_i,
\tilde{\beta}_i'$.
\begin{eqnarray*}
&& g_\theta =c_\theta =I_n,\quad
(-g_\theta\otimes g_\theta, \bc_\theta\otimes c_\theta)
= (- I_{n^2}, I_{n^2}).\\
&& \cM_0(\theta) = \LL_{-1,1}^{\oplus n^2},\quad
H^0(C_0,\cM_0(\theta)) = 0 = H^1(C_0,\cM_0(\theta)).
\end{eqnarray*}
\item[Case 2.] $\lambda=\tilde{\delta}$.
\begin{eqnarray*}
&& g_\theta = \ba_\theta,\quad c_\theta=I_n,\quad
(-g_\theta \otimes g_\theta, \bc_\theta\otimes c_\theta) =
(-\ba_\theta\otimes \ba_\theta, I_{n^2}).\\
&& \cM_0(\theta)= \LL_{-e^{-2\sqrt{-1}\theta},1} \oplus
\LL_{-e^{-\sqrt{-1}\theta},1}^{\oplus 2(n-1)} \oplus
\LL_{-1,1}^{\oplus (n-1)^2}.
\end{eqnarray*}

\item[Case 3.] $\lambda=\tilde{\gamma}$.
\begin{eqnarray*}
&& g_\theta = I_n, \quad c_\theta = a_\theta,\quad
(-g_\theta\otimes g_\theta, \bc_\theta\otimes c_\theta)
=(-I_{n^2}, \ba_\theta\otimes a_\theta).\\
&& \cM_0(\theta)= \LL_{-1,e^{-\sqrt{-1}\theta}}^{\oplus (n-1)}
\oplus
\LL_{-1,e^{\sqrt{-1}\theta}}^{\oplus (n-1)} \oplus \LL_{-1,1}^{\oplus(n^2-2n+2)}.\\
&& H^0(C_0, \cM_0(\theta))= 0 = H^1(C_0, \cM_0(\theta)).
\end{eqnarray*}
\end{enumerate}

In Case 1 and Case 3, \eqref{eqn:K-long} is reduced to
\begin{equation}\label{eqn:K-short}
0 \to \bC_{p_+}^{n^2} \oplus \bC_{p_-}^{n^2}\to
H^1(C,\cM_{\theta,0})\to H^1(C_+,\cM_+(\theta))\oplus
H^1(C_-,\cM_-(\theta))\to 0
\end{equation}
The involution $\hat{\tau}_{\theta,0}$ acts on the exact sequence
\eqref{eqn:K-short} in the following way :
\begin{eqnarray*}
&& \hat{\tau}_{\theta,0}: \bC_{p_+}^{n^2} \oplus \bC_{p_-}^{n^2} \to
\bC_{p_+}^{n^2}\oplus \bC_{p_-}^{n^2},
\quad  (v, w)\mapsto (\bar{w}, \bar{v}) \quad\quad (\mbox{independent of } \theta)\\
&& H^1(C,\cM_{\theta,0})
\stackrel{\hat{\tau}_{\theta,0}}{\longleftrightarrow}
H^1(C,\cM_{\theta,0}),\quad
H^1(C_+,\cM_+(\theta))\stackrel{\hat{\tau}_{\theta,0}}{\longleftrightarrow}
H^1(C_-,\cM_-(\theta)).
\end{eqnarray*}
Thus, we have
$$
0\to \bC^{\oplus n^2} \to
H^1(C,\cM_{\theta,0})^{\hat{\tau}_{\theta,0}} \to
H^1(C_+,\cM_+(\theta))\to 0.
$$
Therefore $i_0^*E \cong E_\bC$, where $E_\bC\to S^1$ is a complex
vector bundle whose fiber at $\theta\in S^1$ is
$H^1(C_+,\cM_+(\theta)) \oplus \bC^{\oplus n^2}$. Therefore $i_0^*E$
is orientable.

\medskip

In Case 2, \eqref{eqn:K-long} is reduced to
\begin{equation}\label{eqn:K-median}
\begin{aligned}
0 \to& H^0(C_0,\cM_0(\theta))
\to \bC_{p_+}^{n^2}\oplus \bC_{p_-}^{n^2}\to H^1(C,\cM_{\theta,0})\\
\to& H^1(C_+,\cM_+(\theta))\oplus H^1(C_0,\cM_0(\theta))\oplus
H^1(C_-,\cM_-(\theta))\to 0
\end{aligned}
\end{equation}
Taking fixed points of the involution  $\hat{\tau}_{\theta,0}$  on
\eqref{eqn:K-median} yields
\begin{eqnarray}\label{eqn:Kreal}
0 \to H^0(C_0,\cM_0(\theta))^{\hat{\tau}_{\theta,0}} \to
\bC_{p_+}^{n^2} \to H^1(C,\cM_{\theta,0})^{\hat{\tau}_{\theta,0}}\to\notag\\
H^1(C_+,\cM_+(\theta))\oplus
H^1(C_0,\cM_0(\theta))^{\hat{\tau}_{\theta,0}}\to 0
\end{eqnarray}
where $\bC_{p_+}^{n^2}$ and $H^1(C_+,\cM_+(\theta))$ are complex
vector spaces.

Recall from \eqref{eqn:lambda-zero} that the holonomies of
$\cM_0(\theta)$ along $\delta$, $\gamma$, $\delta'$, $\gamma'$ are
given by
$$
(-a_\theta\otimes a_\theta, T, -\ba_\theta\otimes \ba_\theta,T).
$$
Let $\tilde{\delta}_0: S^1 \to
Z_{\mathrm{flat}}^{0,2}(U(n^2))^{\hat{\tau}}$ be the loop defined by
\begin{equation}\label{eqn:delta-zero}
\tilde{\delta}_0(\theta) = (-a_\theta\otimes a_\theta, T,
-\ba_\theta\otimes \ba_\theta,T),
\end{equation}
then the loop $\tilde{\delta}_0$ can be viewed as the $S^1$-family
of the bundles $\cM_0(\theta)$. Let $U_n^\vir$ be defined as in
Definition \ref{def:Uvir}. Then the fiber of (the pull back bundle)
$E^{\vir}_0 = \tilde{\delta}_0^*
\left((U_n^\vir)^{\hat{\tau}}\right)$ at $\theta\in S^1$ is the
virtual vector space
$$
H^1(C_0,\cM_0(\theta))^{\hat{\tau}_{\theta,0} } -H^0(C_0,\cM_0(\theta))^{\hat{\tau}_{\theta,0} }.
$$
To show that $\tilde{\delta}^*W_{n,k}\to S^1$ is orientable,
it remains to show that the virtual real vector bundle
$E^\vir_0\to S^1$ is orientable. This is true by Lemma \ref{thm:Czero} below.
\end{proof}

\begin{lm}\label{thm:Czero}
$E^\vir_0$ is an orientable virtual real vector bundle over $S^1$.
\end{lm}

\begin{proof}
We use the notation in the above proof.

\noindent
{\em Step 1.} Recall that $T$ is an involution
on $\bC^n\otimes \bC^n =\bC^{n^2}$ defined by
$T(u\otimes v) = v\otimes u$ for
$u,v\in \bC^n$. The eigenspaces of $T$ are
$E_{+1} =\mathrm{Sym}^2(\bC^n)$ and $E_{-1}=\Lambda^2(\bC^n)$.
We have an orthogonal decomposition
$$
\bC^n\otimes \bC^n = E_{+1} \oplus E_{-1}.
$$
Let $I_{E_{\pm 1}}:E_{\pm 1}\to E_{\pm 1}$ be the identity map.
Then
$$
I_{n^2} = I_{E_{+1}}\oplus I_{E_{-1}},\quad
T= I_{E_{+ 1}} \oplus (-I_{E_{-1}}).
$$
Define $\chi: [0,1]\to U(n^2)$ by
$$
\chi(t) = I_{E_{+1}}\oplus (-e^{i\pi t} I_{E_{-1}}).
$$
Then
$$
\chi(0)= T,\quad \chi(1)= I_{n^2},\quad
\overline{\chi(t)} =\chi(t)^{-1}=\chi(-t).
$$

Note that $a_\theta\otimes a_\theta$ is of the form
$a_\theta\otimes a_\theta =A^+_\theta  \oplus A^-_\theta$,
where $A^+_\theta$ and $A^-_\theta$ are linear automorphisms on $E_{+1}$
and $E_{-1}$, respectively. Therefore
$a_\theta\otimes a_\theta$ and $\chi(t)$ commute
for any $\theta\in S^1$ and $t\in [0,1]$.

\noindent
{\em Step 2.}
For any $t\in [0,1]$ and $\theta\in S^1$, define
\begin{equation}\label{eqn:delta-t}
\tilde{\delta}_t(\theta) =
(-a_\theta\otimes a_\theta, \chi(t), -\ba_\theta\otimes\ba_\theta, \overline{\chi(t)}).
\end{equation}
In particular,  $\tilde{\delta}_0(\theta)$ is given by \eqref{eqn:delta-zero}.
By Step 1, the right hand side of \eqref{eqn:delta-t}
lies in $Z_{\mathrm{flat}}^{0,2}(U(n^2))^{\hat{\tau}}$. So
for each $t\in [0,1]$, \eqref{eqn:delta-t} defines a loop
$\tilde{\delta}_t: S^1 \to Z_{\mathrm{flat}}^{0,2}(U(n^2))^{\hat{\tau}}$.
Let $E^\vir_1 = \tilde{\delta}_1^*
\left((U_n^\vir)^{\hat{\tau}}\right)\to S^1$.
The loops $\tilde{\delta}_0$ and
$\tilde{\delta}_1$ are homotopic, so $E^\vir_0$ is
orientable if and only if $E^\vir_1$ is orientable:

\noindent
{\em Step 3.} Let $\D(\theta)$ be the rank $n^2$, degree
$0$ polystable vector bundle over $C_0$ whose  holonomies along
$\delta, \gamma,\delta', \gamma'$ are given by
\begin{equation}\label{eqn:delta-one}
\tilde{\delta}_1(\theta) =
(-a_\theta\otimes a_\theta, I_{n^2}, -\ba_\theta\otimes\ba_\theta,
I_{n^2}).
\end{equation}
The holonomies of $\D(\theta)$ along $\gamma,
\gamma'$ are $I_{n^2}$, so we may degenerate the torus $C_0$ to a
rational nodal curve $C'$ (see Figure 3), and degenerate
$\D(\theta)$ to $\D'(\theta)$.
\begin{figure}[h]
\begin{center}
\psfrag{p+}{\tiny $p_+$}
\psfrag{p-}{\tiny $p_-$}
\psfrag{p}{\tiny $p$}
\psfrag{r}{\tiny $\gamma$}
\psfrag{rp}{\tiny $\gamma'$}
\psfrag{d}{\tiny $\delta$}
\psfrag{dp}{\tiny $\delta'$}
\psfrag{nu}{\tiny $\nu$}
\psfrag{p0}{\tiny $q_+$}
\psfrag{p1}{\tiny $q_-$}
\psfrag{C0}{\small $C_0$}
\psfrag{Cp}{\small $C'$}
\includegraphics[scale=0.6]{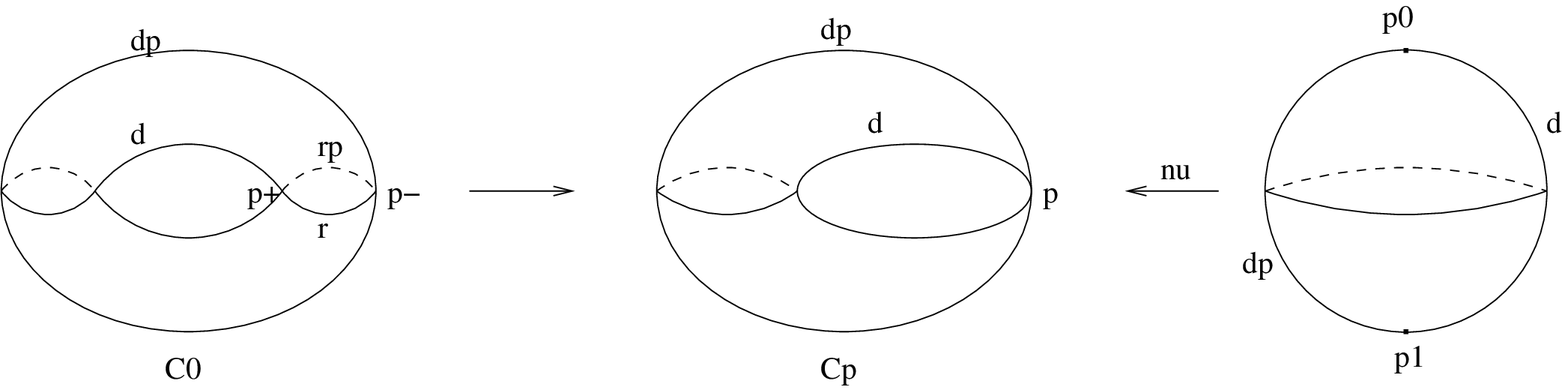}
\end{center}
\caption{Degeneration of $C_0=S^1\times S^1$. }
\end{figure}

More explicitly, we consider  a family of
Riemann surface $R_s$, where $s\in I=[0,1]$, such that
\begin{enumerate}
\item[(i)] $R_1=C_0$, $R_0=C'$.
\item[(ii)] $R_s$ is a smooth torus for $0<s\leq 1$.
\item[(iii)] There is a
family of anti-holomorphic involution $\sigma_s:R_s \to R_s$, such
that $\sigma_0(p)=p$, $\sigma_1=\tau_0$, $\sigma_s(\delta)=\delta'$,
and $\sigma_s(\gamma)=\gamma'$ if $s\neq 0$.
\item[(iv)] There is a normalization map $\nu: \bP^1 \to C'$
such that $\nu(q_\pm)=p$ (see Figure 3).
There is an anti-holomorphic involution $\tilde{\sigma}:\bP^1\to \bP^1$
such that
$$
\nu\circ \tilde{\sigma}=\sigma_0\circ \nu,\quad
\tilde{\sigma}(q_\pm)=q_\mp.
$$
\end{enumerate}
In the degeneration $s\to 0$, the loop $\gamma\gamma'$
shrinks to a point $p$.

We consider a family of  polystable holomorphic
vector bundles $\D_{\theta,s}$ of rank $n^2$, degree $0$
over $R_s$ with the following properties:
\begin{enumerate}
\item[(i)] For $0< s\leq 1$, the holonomies of $\D_{\theta,s}$
along $\delta, \gamma,\delta',\gamma'$ are given by
\eqref{eqn:delta-one},
so we have $\overline{\sigma_s^* \D_{\theta,s}}=\D_{\theta,s}$.
\item[(ii)] When $s=0$, we have $\D_{\theta,0} =\D'(\theta)$, and
$\nu^* \D'(\theta) =\cO_{\bP^1}^{\oplus n^2}$.
$\D'(\theta)$ is obtained by identifying $v\in \bC^{n^2}_{q_+}$
with $(-a_\theta\otimes a_\theta) v \in  \bC^{n^2}_{q_-}$.
\item[(iii)] For all $(\theta,s)\in S^1\times I$, there
is a conjugate linear involution $\sigma_{\theta,s}$ on
$H^*(R_s, \D_{\theta,s})$ such that
$$
H^1(R_1, \D_{\theta,1})^{\sigma_{\theta,1}}-
H^0(R_1,\D_{\theta,1})^{\sigma_{\theta,1}}
$$is the fiber of $E_1^\vir$ at $\theta\in S^1$.

\item[(iv)] The conjugate linear involution $\sigma_{\theta,0}$ on $H^0(R_0,\D_{\theta,0})$ induces
a conjugate linear involution $\sigma_{\theta,0}$ on
$H^0(\bP^1,\cO_{\bP^1}^{\oplus n^2})\cong \bC^2$ (constant sections)
that is given by $v\mapsto \bar{v}$.
\end{enumerate}

As $(\theta,s)$ varies, the virtual real vector spaces
$$
\tilde{E}_{\theta,s}^{\vir}= H^1(R_s,\D_{\theta,s})^{\sigma_{\theta,s}}
-H^0(R_s,\D_{\theta,s})^{\sigma_{\theta,s}}
$$
form a virtual real vector bundle $\tilde{E}^\vir$ over
the cylinder $S^1\times I$. Let $i_s: S^1\to S^1\times I$
be the embedding $\theta \mapsto (\theta,s)$. Then
$$
i_0^* \tilde{E}^\vir = F^\vir,\quad i_1^* \tilde{E}^\vir= E_1^\vir.
$$
where the fiber of $F^\vir$ at $\theta\in S^1$ is
\begin{equation}\label{eqn:F}
(F^\vir)_\theta= H^1(C', \D'(\theta))^{\sigma_{\theta,0}}
-H^0(C',\D'(\theta))^{\sigma_{\theta,0}}.
\end{equation}
The maps $i_0$ and $i_1$ are homotopic, so
$E^\vir_1$ is orientable if and only if $F^\vir$
is orientable.

\noindent
{\em Step 4.} We have a long exact sequence
\begin{equation}\label{eqn:Pone}
0 \to H^0(C',\D'(\theta)) \to H^0(\bP^1,\cO_{\bP^1}^{\oplus n^2} )
\cong \bC^{n^2} \stackrel{f_\theta}{\to} \bC^{n^2}_p
  \to H^1(C',\D'(\theta)) \to  0
\end{equation}
where $f_\theta(v) =  v + (a_\theta \otimes a_\theta)v$ and we used
$H^1(\bP^1,\cO_{\bP^1}^{\oplus n^2} ) =0$. The involution
$\sigma_{\theta,0}$ acts on \eqref{eqn:Pone} by
\begin{eqnarray*}
&H^0(\bP^1,\cO_{\bP^1}^{\oplus n^2})\cong\bC^{n^2}\to
H^0(\bP^1,\cO_{\bP^1}^{\oplus n^2})& ,\quad
v\mapsto \bar{v},\\
&\bC^{n^2}_p\to \bC^{n^2}_p&,\quad v\mapsto (a_\theta\otimes a_\theta)\bar{v}.
\end{eqnarray*}
Therefore
\begin{eqnarray*}
H^0(\bP^1,\cO_{\bP^1}^{\oplus n^2})^{\sigma_{\theta,0}}  &=& \bR^{n^2}\\
\left(\bC^{n^2}_p\right)^{\sigma_{\theta,0}} &=&
e^{\sqrt{-1}\theta} \bR \oplus
e^{\sqrt{-1}\theta/2}\bR^{2(n-1)}
\oplus \bR^{(n-1)^2}
\end{eqnarray*}
So we have an exact sequence
\begin{equation} \label{eqn:ER}
0\to H^0(C',\D'(\theta))^{\sigma_{\theta,0}}\to\bR^{n^2}
\stackrel{f_\theta}{\to} E_{\bR,\theta} \to
H^1(C',\D'(\theta))^{\sigma_{\theta,0}}\to 0.
\end{equation}
where
$$
E_{\bR,\theta} = e^{\sqrt{-1}\theta} \bR \oplus
e^{\sqrt{-1}\theta/2}\bR^{2(n-1)} \oplus \bR^{(n-1)^2}.
$$

To show that $F^\vir$ is orientable, it suffices to show that
$E_{\bR,\theta}$ form an orientable real vector bundle $E_\bR \to
S^1$. Let $L_0$ and $L_1$ denote the trivial and nontrivial real
line bundles over $S^1$, so that
$$
w_1(L_j)= j \in \bZ/2\bZ=H^1(S^1,\bZ/2\bZ),\quad j=0,1.
$$
Then
$$
E_\bR = L_0^{\oplus(n^2-2n+2)} \oplus L_1^{\oplus(2n-2)}, \quad
w_1(E_\bR)=0 \in \bZ/2\bZ = H^1(S^1,\bZ/2\bZ).
$$
Therefore $E_\bR\to S^1$ is orientable.

\end{proof}

\section*{Acknowledgments}
The first author was partially supported by NSC 97-2628-M-006-013-MY2.
The second author was partially supported by the Sloan Research Fellowship.
The third author was partially supported by NSF grants DMS-0343640
(RTG), DMS-0804553, and DMS-0968766.

\noindent
{\sc Nan-Kuo Ho\\
Department of Mathematics\\
National Tsing Hua University and\\
National Center for Theoretical Sciences\\
Hsinchu 300\\
Taiwan} \\
{\em E-mail address:} {\tt nankuo@math.nthu.edu.tw}

\bigskip

\noindent
{\sc Chiu-Chu Melissa Liu\\
Department of Mathematics\\
Columbia University\\
New York\\
New York 10027 \\
USA}\\
{\em E-mail address:} {\tt ccliu@math.columbia.edu}

\bigskip

\noindent
{\sc Daniel Ramras\\
Department of Mathematical Sciences\\
New Mexico State University\\
Las Cruces\\
New Mexico 88003-8001\\
USA}\\
{\em E-mail address:} {\tt ramras@nmsu.edu}

\end{document}